\documentclass{article}

\usepackage[top=4cm, bottom=3cm, left=2.5cm, right=2.5cm]{geometry}
\usepackage[utf8]{inputenc}  
\usepackage[T1]{fontenc}     
\usepackage{hyperref}        
\usepackage{url}             
\usepackage{booktabs}        
\usepackage{amsfonts}        
\usepackage{nicefrac}        
\usepackage{microtype}       
\usepackage{lipsum}
\usepackage{graphicx}
\graphicspath{ {./images/} }
\usepackage{amssymb}
\usepackage{latexsym}
\usepackage{graphicx}  
\usepackage{subcaption}  
\usepackage{amsmath}  
\usepackage{array}
\usepackage{hhline}
\usepackage{url}
\usepackage{xcolor}
\definecolor{newcolor}{rgb}{.8,.349,.1}

\usepackage{framed,multirow}
\usepackage{stmaryrd}

\usepackage{algorithm}
\usepackage{algpseudocode}
\usepackage{amsmath}

\usepackage[numbers, sort&compress]{natbib} 
\usepackage{amsmath}
\numberwithin{equation}{section}

\title{ The modified Physics-Informed Hybrid
Parallel Kolmogorov--Arnold and Multilayer
Perceptron Architecture with domain decomposition 
}

\author{
	Qiumei Huang$^{\,1}$, Xu Wang$^{\,1}$, Yu Zhao$^{\,1, *}$ }
\date{}

\begin{document}
	 \maketitle
	\begin{center}
		{\footnotesize 1 School of Mathematics, Statistics and Mechanics, Beijing University of Technology, Beijing, 100124, China. \\
		}
	\end{center}
\begingroup
\renewcommand{\thefootnote}{} 
\footnotetext{
* Corresponding author.\\
\hspace*{2.5em}E-mail: qmhuang@bjut.edu.cn (Qiumei Huang), wx202366070@emails.bjut.edu.cn (Xu Wang), y.zhao@bjut.edu.cn (Yu Zhao).
}
\endgroup

\begin{abstract}
In this work, we propose a modified Hybrid Parallel Kolmogorov--Arnold Network and Multilayer Perceptron Physics-Informed Neural Network to overcome the high-frequency and multiscale challenges inherent in Physics-Informed Neural Networks. This proposed model features a trainable weighting parameter to optimize the convex combination of outputs from the Kolmogorov--Arnold Network and the Multilayer Perceptron, thus maximizing the networks' capabilities to capture different frequency components. Furthermore, we adopt an overlapping domain decomposition technique to decompose complex problems into subproblems, which alleviates the challenge of global optimization. Benchmark results demonstrate that our method reduces training costs and improves computational efficiency compared with manual hyperparameter tuning in solving high-frequency multiscale problems.

\end{abstract}


\section{ Introduction }

Traditional numerical methods for solving partial differential equations, such as finite volume method and finite element method, are largely based on mesh discretization. However, these methods often encounter numerous difficulties when dealing with complex geometries or high-dimensional problems. To overcome these limitations, researchers have explored several mesh-free methods, along with neural network based methods have emerged as a promising direction. In particular, Physics-Informed Neural Networks (PINNs) have gained extensive attention \cite{1}. PINNs employ neural networks to construct a continuous parametric representation of the solution, and are trained by minimizing loss functions
that incorporate the physics constrains and the initial and boundary conditions. This framework not only ensures high accuracy of the solution but also eliminates the complexity associated with mesh generation in complex geometrical domains \cite{5}.
Moreover, PINNs have an advantage in solving continuous governing equations by avoiding the need for mesh discretization and reducing the mathematical constraints on boundary conditions, which is effective in handling problems with irregular geometries and insufficient constraints.
Based on these benefits, PINNs have been applied to solve partial differential equations in various fields, including fluid dynamics, cardiac electrophysiology, wave propagation, and thermal conduction \cite{16,17,18}.

Despite the success in many applications, PINNs still encounter several fundamental limitations in solving partial differential equations. Firstly, PINNs predominantly employ the Multilayer
Perceptron (MLP) as the base architecture. The use of MLP restricts PINNs' ability to capture nonlinear features in complex physical domains, thus affecting their accuracy and generalization and limiting their applications. 
Furthermore, although recent studies have begun to bridge fundamental theoretical gaps, the convergence properties lack rigorous theoretical analysis compared to classical numerical methods \cite{21}.
More importantly, PINNs suffer from the spectral bias inherent in neural networks, which results in lower learning efficiency for high-frequency components than low-frequency components \cite{23}.
This bias leads to a gap in convergence rates \cite{25,26}, which becomes particularly pronounced when solving problems with high-frequency and multiscale features \cite{27,28}, and can often result in a sharp drop in the solution accuracy of PINNs. To have an acceptable accuracy, PINNs need to increase the number of collocation points and expand the neural network size, which results in a superlinear growth in the number of degrees of freedom and a corresponding surge in computational cost \cite{29}. Thus, traditional numerical methods generally outperform PINNs in forward modeling, and their current applications are largely limited to small-scale problems characterized by low-frequency features.

To improve the performance of PINNs in complex physical problems, researchers have proposed various improved algorithms for the PINNs. One of the key focuses is to adopt a neural network architecture with enhanced representational capabilities. Among these networks, Convolutional Neural Networks (CNNs) \cite{41} and Graph Neural Networks (GNNs) \cite{42} have been explored as alternatives to traditional MLP. Furthermore, an important advancement is the Kolmogorov--Arnold Network (KAN) proposed by Li et al. \cite{58}. Based on the Kolmogorov--Arnold representation theorem, KAN expresses an arbitrary multivariate continuous function as a finite composition of univariate functions. The hierarchical function composition mechanism of KAN aligns intrinsically with the superposition principle of differential equation solutions. The mechanism not only improves the ability to approximate highly nonlinear problems but also enables seamless integration of physical priors and data-driven learning through a combination of learnable basis functions.
Consequently, the adoption of KAN enhances model interpretability and expressive power for complex physical phenomena, providing a more robust and reliable framework for modeling intricate physical systems. In addition to the network architecture, another significant focus is to enhance the optimization process of PINNs, such as adaptive weighting methods \cite{45}, residual-based attention mechanisms \cite{46}, adaptive weights method based on balanced residual decay rate \cite{48}, stacked networks \cite{stacked}, causal constraint methods \cite{49}, multi-fidelity PINNs \cite{51}, progressive hierarchical learning frameworks \cite{53}, Fourier Feature Networks \cite{30}, Multi-scale Deep Neural Networks \cite{32} and Pre-trained Physics-Informed Neural Networks (PT-PINN) \cite{zengjia33PTPINN}. These methods improve the performance and convergence by introducing additional mechanisms (such as adaptation, attention, sampling, etc.), but they increase the complexity of the algorithms and the reliance on hyperparameters or prior knowledge.

As another improvement to PINNs, integrated domain decomposition technique has been successfully applied to solve problems involving multi-scale and high-frequency phenomena. PINNs with domain decomposition technique divide complex computational domains into multiple smaller subdomains for distributed solutions, which significantly enhances the computational speed, robustness, and solution accuracy of various deep learning models. Thus, the strategy of integrating various improved PINNs with domain decomposition techniques has become a computational paradigm that attracts increasing attention due to its superior convergence rate and higher computational accuracy.
For example, the extended Physics-Informed Neural Network (XPINN) \cite{40} employs a non-overlapping domain decomposition strategy, where independent neural networks are trained on separate subdomains. However, the absence of overlapping regions between subdomains introduces discontinuities at subdomain boundaries, necessitating additional coupling loss terms to facilitate information exchange between subdomains. 
In contrast, the Finite Basis Physics-Informed Neural Networks (FBPINNs) \cite{39} \textcolor{black}{use partition of unity functions to weight the output of neural networks in each subdomain,} eliminating coupling loss terms and ensuring solution continuity across subdomains for better accuracy. Building upon FBPINNs, the stacking FBPINN approach \cite{timedd} integrates multifidelity stacking PINNs with time domain decomposition FBPINNs to address time-dependent multiscale problems. Furthermore, multilevel FBPINNs \cite{MDDPINN}, by incorporating multilevel overlapping domain decomposition into FBPINNs, significantly enhance the accuracy and computational efficiency in solving high-frequency and multiscale problems.

In this work, we propose a modified Hybrid Parallel Kolmogorov--Arnold Network and Multilayer Perceptron Physics-Informed Neural Network (modified HPKM-PINN) model based on overlapping domain decomposition to solve high-frequency oscillatory multiscale partial differential equations. To combine MLP's strength in capturing low-frequency features with KAN's capacity for resolving high-frequency details, we construct a convex combination of MLP and KAN outputs to balance their contributions in handling different frequency features and reduce computational costs. This combination is designed by using trainable weight parameters, which are optimized through S-shaped functions to ensure strict confinement within the interval [0,1], thereby enabling adaptive modulation of the contribution ratio between MLP and KAN outputs. The adaptive convex combination achieves optimization without manual tuning and reduces computational costs.
Furthermore, the overlapping domain decomposition technique decomposes complex problems into subproblems to overcome the difficulties in global optimization. A variety of benchmark results demonstrate that our modified model reduces training costs and improves computational efficiency when solving high-frequency multiscale problems. The proposed modified model establishes a novel computational framework for solving high-frequency partial differential equations.

The rest of this paper is organized as follows. Section 2 reviews the Hybrid Parallel Kolmogorov--Arnold and MLP (HPKM) model and presents our modified HPKM-PINN framework integrating overlapping domain decomposition. Section 3 presents extensive numerical experiments to investigate the network architecture and hyperparameter selection strategy, and provide a comprehensive analysis of the modified HPKM-PINN model's performance under various parameter configurations for solving Poisson, Helmholtz, reaction-diffusion, and Allen--Cahn equations. Finally, the conclusion and future work are discussed in Section 4.

\section{Method}

The theoretical framework is developed in two steps: in subsection 2.1, by introducing an adaptively adjustable parameter, we propose a modified HPKM architecture to automatically balance the contributions of KAN's and MLP's inputs, thereby adaptively dealing with variety of frequency in the problem; in subsection 2.2, we integrate this architecture with the overlapping domain decomposition technique to decompose the global problem into various local problems and to compute parallelly, which improves the computational efficiency and solution accuracy.

\subsection{The modified HPKM architecture}

This subsection first introduces the two core components of the HPKM architecture, namely the MLP and KAN. A modified HPKM architecture is then proposed based on these components.

\subsubsection{MLP and KAN}

In PINNs, the predictive solution is typically constructed using an MLP architecture. As a feedforward neural network, the MLP framework consists of fully connected layers arranged in a cascaded structure, including the input, hidden, and output layers. It is designed to perform nonlinear feature extraction and spatial mapping of the input data.

The parameters \(\theta\) of an MLP architecture are denoted as \(\theta=\{W^{(i)}, b^{(i)}\}_{i=1}^n\), where $n$ denotes the number of layers, \(W^{(i)}\) and \(b^{(i)}\) represent the weight matrix and bias of the $i$-th layer, respectively. The output \( u(\mathbf{x};\theta) \) of the MLP with $n$ layers can be expressed as
\begin{equation*}
u(\mathbf{x};\theta) = \sigma \left( W^{(n)} \sigma \left( W^{(n-1)} \cdots \sigma \left( W^{(1)} \mathbf{x} + b^{(1)} \right) \cdots + b^{(n-1)} \right) + b^{(n)} \right),
\end{equation*}
where \(\sigma\) is the activation function, $\mathbf{x} = (x_1, \dots,x_d)^T \in \mathbb{R}^d$ is the input, and $d$ is the input dimension. According to the Universal Approximation Theorem \cite{57}, such an architecture with adequate neurons and appropriate activation functions can
approximate any continuous function defined on compact subsets of $\mathbb{R}^d$.

Based on the Kolmogorov--Arnold Representation Theorem, Liu et al. proposed an approximate model of the multivariate function $f(\mathbf{x})$ in the following form \cite{58} 
\begin{equation*}
    f(\mathbf{x}) \approx \sum_{i_{L-1}=1}^{m_{L-1}} \varphi_{L-1, i_{L}, i_{L-1}} \left( \sum_{i_{L-2}=1}^{m_{L-2}} \cdots \left( \sum_{i_2=1}^{m_2} \varphi_{2, i_3, i_2} \left( \sum_{i_1=1}^{m_1} \varphi_{1, i_2, i_1} \left( \sum_{i_0=1}^{m_0} \varphi_{0, i_1, i_0} \left( x_{i_0} \right) \right) \right) \right) \cdots \right),
\end{equation*}
where $\varphi_{i,j,k}$ is the univariant activation function in the $i$-th layer that connects the $j$-th input node to the $k$-th output node, $L$ is the number of layers in the KAN, and $\{m_i\}_{i=0}^{L}$ denotes the number of nodes in the $i$-th layer. The performance variations among KAN architectures are governed by the specific selection of the basis functions $\varphi_{i,j,k}(x)$.

The univariate activation function in the original KAN is formulated as a weighted sum of the basis function $b(x)$ and the B-spline function $\text{spline}(x)$
\begin{equation*}
\varphi(x) = w_b b(x) + w_s \text{spline}(x), \quad\quad~ \text{with}~~~~b(x) = \frac{x}{1 + e^{-x}},\quad
\text{spline}(x) = \sum_{i} c_i B_i(x),
\end{equation*}
where $B_i(x)$ is a polynomial of degree $k$. \( \text{spline}(x) \) denotes the \( k \)-th order B-spline function defined on a uniform grid partitioned into \( G \) intervals over the domain. For the one-dimensional domain \([r_{min}, r_{max}]\), the grid with \( G \) intervals has points \( \{t_0 = r_{min}, t_1, t_2, \dots, t_G = r_{max}\} \); see \cite{58}. $w_b$
, $w_s$, and $c_i$ are trainable parameters. Through the dynamic adjustment of the parameters $\theta = \{ w_b, w_s, c_i \}$, KAN acquires the capability to learn specific adaptive activation functions instead of the predefined static nonlinearity of MLP.

Recently, many KAN architectural variants have been proposed. In particular, the Radial Basis Function (RBF) KANs \cite{55} and Chebyshev KANs \cite{56} maintain compatibility with domain decomposition technique while achieving superior computational efficiency and approximation accuracy compared to conventional B-spline implementations. Based on these developments, we employ a Fourier KAN architecture \cite{FKAN} that replaces B-spline basis with Fourier series, achieving enhanced performance in high-frequency component approximation compared to conventional approaches. The Fourier series are denoted as follows 
\begin{equation*}
  \textcolor{black}{  \psi(x) = \frac{a_0}{2} + \sum_{j=1}^K (a_j \cos(jx) + b_j\sin(jx)) },
\end{equation*}
where $K$ is the number of frequency components, \(a_j, j = 0, \cdots, K\) and \(b_j, j = 1, \cdots, K\) denote the learnable Fourier coefficients.

\subsubsection{ Modified HPKM architecture}

By employing a parallel structure, the HPKM architecture is proposed to integrate KAN's high-frequency detail resolution capability with MLP's low-frequency global approximation capacity \cite{zengjiaHPKM}. The output of HPKM is formulated as a weighted combination of the two networks' outputs,
\begin{equation}
\label{method: alpha}
  u(\mathbf{x};\theta) = \alpha ~u_{\text{KAN}}(\mathbf{x};\theta_{\text{KAN}}) + (1-\alpha)~u_{\text{MLP}}(\mathbf{x};\theta_{\text{MLP}}),
\end{equation}
where $\mathbf{x}$ is the input vector, \(\alpha\) is a pre-defined tunable weight parameter, $\theta = \{\theta_{\text{KAN}}, \theta_{\text{MLP}}\}$ is the set of total parameters of the HPKM model with $\theta_{\text{KAN}}$, $ \theta_{\text{MLP}}$ the parameters of KAN and MLP, respectively. $u_{\text{KAN}}(\mathbf{x};\theta_{\text{KAN}})$ and $u_{\text{MLP}}(\mathbf{x};\theta_{\text{MLP}})$ in the output expression \eqref{method: alpha} are the outputs of KAN and MLP, respectively.

The value of \(\alpha\) requires manual testing and selection for optimal performance across different problems. To overcome this limitation, we reformulate \(\alpha\) as a trainable parameter that dynamically adjusts during optimization.
This modification eliminates the need for laborious hyperparameter tuning, while enabling the adaptive balancing of relative contributions between MLP and KAN components throughout the training process. Meanwhile, we employ an S-shaped function \(S(\alpha)\) to bound the weighting parameters to the interval \([0,1]\), producing a convex combination of outputs and enabling more stable optimization. Thus, one gets the output of the modified HPKM architecture
\begin{equation}
\label{method: Salpha}
  u_{\text{MHPKM}}(\mathbf{x};\theta) = S(\alpha) u_{\text{KAN}}(\mathbf{x};\theta_{\text{KAN}}) + (1-S(\alpha))u_{\text{MLP}}(\mathbf{x};\theta_{\text{MLP}}). 
\end{equation}

\begin{figure}[htbp]     
    \centering 
    \includegraphics[scale=0.7]{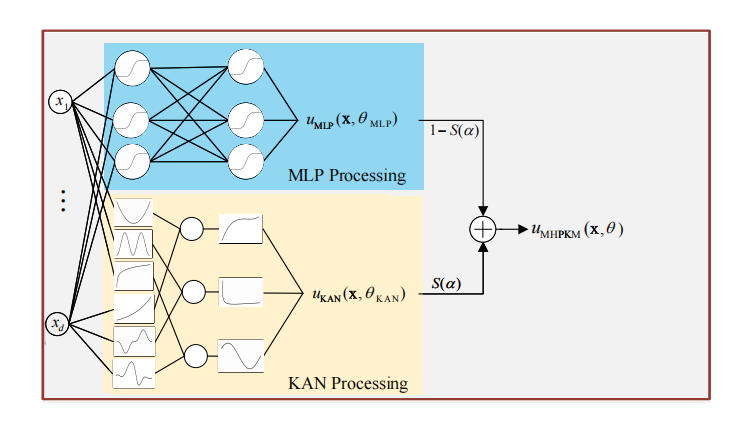}
    \caption{Modified HPKM architecture.} \label{figure1}
\end{figure}

The adaptive weighting coefficient $S(\alpha)$ in \eqref{method: Salpha} provides the following advantages: (i) automatically balancing MLP's global approximation capability with KAN's local feature extraction for multi-frequency modeling; (ii) eliminating manual parameter tuning through adaptive \(\alpha\) adjustment while reducing computational costs; and (iii) guaranteeing stable training dynamics through the constrains of the S-shaped function \(S(\alpha) \in [0,1] \). Figure $\ref{figure1}$ shows the design of the modified HPKM architecture.  

\subsection{Overlapping domain decomposition technique in modified HPKM-PINN}

In this subsection, we apply the modified HPKM architecture to the PINN framework, introducing a modified HPKM-PINN that incorporates the overlapping domain decomposition technique.

\subsubsection{Overlapping domain decomposition}
\label{domaindecomposition}
The overlapping domain decomposition technique has advantages in handling multiscale and high-frequency problems. In this work, we deploy a neural network on each overlapping subdomain and normalize its inputs. This normalization transforms the high-frequency components within each subdomain into lower-frequency ones, thereby reducing the impact of spectral bias.

\begin{figure}[htbp]     
    \centering 
    \includegraphics[scale=0.6]{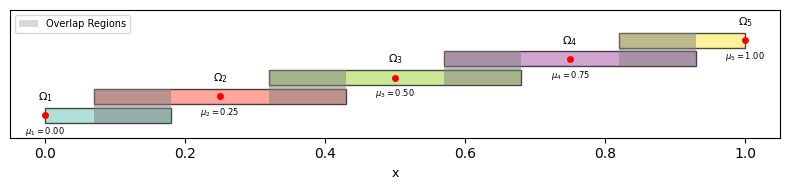}
    \caption{Domain decomposition of the 1D domain $\Omega = [0,l], l=1$.} \label{DD}
\end{figure}

Without loss of generality, we demonstrate the domain decomposition technique on the domain \(\Omega=[0,l_{1}] \times [0,l_{2}] \times \dots\times [0,l_{d}]\), where $l_{j} (j=1,\dots,d$) represents the length in the $j$-th direction of $\Omega$. 
Following the framework established in \cite{39,MDDPINN,FBKAN}, we partition $\Omega$ into $N$ overlapping subdomains \(\{\Omega_i\}_{i=1}^N\), where each \(\Omega_i\) maintains non-empty intersection with its adjacent subdomains, that is 
\begin{equation*}
    \Omega = \bigcup_{i=1}^{N} \Omega_i,~~~~ \Omega_i = \gamma_{i1}\times\gamma_{i2}\times\cdots\times\gamma_{id},
\end{equation*}
where $\gamma_{ij}$ denotes the $x_j$-direction component of the $i$-th subdomain.
The interval $\gamma_{ij}$ with the overlap ratio $\delta_{ij}$ can be expressed as 
\begin{align*}
\gamma_{ij} = \left( ~\max(\mu_{ij}-\delta_{ij}\sigma_{ij},0),~~\min(\mu_{ij}+\delta_{ij}\sigma_{ij},l_{j})~ \right),
\end{align*}
where $\mu_{ij}$ and $\sigma_{ij}$ represent the center and half-width of the $i$-th subdomain in the $x_j$-direction. 
The overlap ratio $\delta_{ij}>1$ reflects the extent of overlap between adjacent regions, with a 
$\delta_{ij}$ indicating a larger overlapping area.
For simplicity, we take a one-dimensional case as an example to illustrate how the domain decomposition is performed. Consider $\Omega = [0,l]$, then
\begin{align*}
\Omega_{i} = \left( ~\max(\mu_i-\delta_i\sigma_i,0), ~\min(\mu_i+\delta_i\sigma_i,l)~ \right),
\end{align*}
with the overlap ratio $\delta_i > 1$, the center $\mu_{i} = \frac{l(i-1)}{N-1}$ and the half-width $\sigma_{i} = \frac{\delta_{i} l}{2(N-1)}$. As shown in Figure \ref{DD}, we take $l=1$, $N=5$, $\delta_i = 1.2, i = 1,\dots ,5$.

Following the overlapping domain decomposition process, we train separate neural networks on each subdomain $\Omega_i$ and construct a global solution by blending their outputs with \textcolor{black}{partition of unity} functions. In this work, the \textcolor{black}{partition of unity} function $\omega_i$ is given by
\begin{align*}
        \omega_i(\mathbf{x}) = \frac{\hat{\omega}_i(\mathbf{x})}{\sum_{j=1}^N \hat{\omega}_j(\mathbf{x})}, ~~\mathbf{x} = (x_1, \dots,x_d)^T \in \Omega, ~~i = 1, \dots, N,
\end{align*}
where for $N = 1$,
\begin{align*}
\hat{\omega}_i(\mathbf{x})\equiv1, ~~ \mathbf{x} \in \Omega,
\end{align*}
and for $N > 1$,
\begin{align*}
    \hat{\omega}_i(\mathbf{x}) = 
\begin{cases} 
\prod_{j = 1}^d(1 + \cos(\pi(x_j - \mu_{ij})/\sigma_{ij}))^2, & \mathbf{x}\in \Omega_i,\\
0, & \mathbf{x}\notin \Omega_i.  
\end{cases}
\end{align*}

Here the \textcolor{black}{partition of unity} function $\omega_i$ constrains the neural network within the \( i \)-th local subdomain. We remark that besides the \textcolor{black}{partition of unity} function employed above, there exist other alternative choices; see \cite{window}. Consequently, $\omega_i$ satisfies $\text{supp}(\omega_i) = \overline{\Omega_i}$, and $\sum_{i=1}^{N} \omega_i(\mathbf{x}) \equiv 1, $ $\forall \mathbf{x} \in \Omega$. Then, for \( i \)-th subdomain, we define the neural network function space $\mathcal{F}_i(\mathbf{x};\theta)$ as 
\begin{equation*}
\mathcal{F}_i(\mathbf{x};\theta) = \left\{f_i(\mathbf{x};\theta^{(i)}) \mid \mathbf{x} \in \Omega_i, \theta^{(i)}\in \Theta_i\right\}, ~~\quad i=1, \dots, N, \quad 
\end{equation*}
where $f_i(\mathbf{x};\theta^{(i)})$ is the neural network function in the $i$-th subdomain, $\theta^{(i)}$ represents the parameters set of the corresponding neural network, $\Theta_i$ is the weight space, representing the set of all trainable parameters of the $i$-th subnetwork. 
Then the final output function can be defined as follows
\begin{equation*}
u(\mathbf{x};\theta) = \sum_{i=1}^N\omega_i(\mathbf{x})f_i(\mathbf{x};\theta^{(i)}),\quad~~\forall \mathbf{x}\in\Omega,
\end{equation*}
where $u(\mathbf{x};\theta)$ represents the final global prediction solution. Due to the construction of the \textcolor{black}{partition of unity} functions, the product $\omega_i(\mathbf{x}) f_i(\mathbf{x};\theta^{(i)})$ vanishes outside $\Omega_i$, thus ensuring a properly localized contribution.

\subsubsection{Modified HPKM-PINN with overlapping domain decomposition}

Under the framework of overlapping domain decomposition, we employ the modified HPKM architecture \eqref{method: Salpha}, and the final output function is     
\begin{equation}    \label{equation-over}
\begin{aligned}
u_{\text{MHPKM}}(\mathbf{x};\theta) 
    &=\sum_{i=1}^N \omega_i(\mathbf{x})u_{\text{MHPKM}}^{(i)}(\mathbf{x};\theta^{(i)})
    \\
    &= \sum_{i=1}^N \omega_i(\mathbf{x}) ( S(\alpha^{(i)})u_{\text{KAN}}^{(i)}(\mathbf{x};\theta_{\text{KAN}}^{(i)}) + (1-S(\alpha^{(i)}))u_{\text{MLP}}^{(i)}(\mathbf{x};\theta_{\text{MLP}}^{(i)})), ~~~\mathbf{x}\in\Omega,
\end{aligned}
\end{equation} %
where $u_{\text{MHPKM}}(\mathbf{x};\theta)$ is the output of the modified HPKM-PINN with domain decomposition, 
$u_{\text{MHPKM}}^{(i)}(\mathbf{x}; \theta^{(i)})$ represents the local prediction on the $i$-th subdomain $\Omega_i$ using the modified HPKM architecture \eqref{method: Salpha}, $u_{\text{KAN}}^{(i)}(\mathbf{x};\theta_{\text{KAN}}^{(i)})$ and $u_{\text{MLP}}^{(i)}(\mathbf{x};\theta_{\text{MLP}}^{(i)})$ represent the output of the KAN and MLP components within the modified HPKM framework in $i$-th subdomain $\Omega_i$, respectively. 
The parameters of all the networks in \eqref{equation-over} are denoted as $\theta = \{\theta^{(i)}\}_{i=1}^N=\{\theta_{\text{KAN}}^{(i)},\theta_{\text{MLP}}^{(i)}, \alpha^{(i)}\}_{i=1}^N$, $\alpha^{(i)}$ is the adaptive weighting parameter for each subdomain network $u_{\text{MHPKM}}^{(i)}(\mathbf{x}; \theta^{(i)})$, and $\omega_i(\mathbf{x})$ represents the corresponding \textcolor{black}{partition of unity} function defined in subsection 2.2.1.

To establish our computational framework for solving partial differential equations, we consider the following system, 
\begin{equation*}
    \begin{cases} 
\frac{\partial}{\partial t}u(\mathbf{x},t) + \mathcal{N}[u(\mathbf{x},t)] = 0, & (\mathbf{x},t) \in \Omega \times (0,T], \\ 
u(\mathbf{x},0) = g(\mathbf{x}), & \mathbf{x} \in \Omega, \\ 
\mathcal{B}[u](\mathbf{x},t) = 0, & (\mathbf{x},t) \in \partial\Omega \times [0,T],
    \end{cases}
\end{equation*}
where \( \mathbf{x} = (x_1, \dots,x_d)^T \in \Omega \subset \mathbb{R}^d \) denotes the spatial coordinate in a bounded domain \(\Omega\) with boundary \(\partial\Omega\), \( t \in [0,T] \) represents the temporal coordinate over the solution interval, \( u(\mathbf{x},t)\) is unknown, \( g(\mathbf{x})\) specifies the initial condition, \(\mathcal{N}[\cdot]\) denotes a differential operator, and \(\mathcal{B}[\cdot]\) denotes the boundary condition operator.

To optimize the modified HPKM-PINN \eqref{equation-over} within the overlapping domain decomposition technique, we formulate the training process as a minimization problem whose loss function is given as
\begin{equation}  \label{eq:loss_function}
    \mathcal{L}(\theta) = \lambda_{ic}\mathcal{L}_{ic}(\theta) + \lambda_{bc}\mathcal{L}_{bc}(\theta) + \lambda_{r}\mathcal{L}_{r}(\theta), 
\end{equation}
with
\begin{equation*}
    \begin{aligned}
        \mathcal{L}_{ic}(\theta) &= \frac{1}{N_{ic}} \sum_{j=1}^{N_{ic}} \left( \sum_{i=1}^{N} \omega_i(\mathbf{x}_{ic}^j)u_{\text{MHPKM}}^{(i)}(\mathbf{x}_{ic}^j; \theta^{(i)}) - g(\mathbf{x}_{ic}^j) \right)^2, \\
        \mathcal{L}_{bc}(\theta) &= \frac{1}{N_{bc}} \sum_{j=1}^{N_{bc}} \left(\mathcal{B}\left[ \sum_{i=1}^{N} \omega_i(\mathbf{x}_{bc}^j  ) u_{\text{MHPKM}}^{(i)}(\mathbf{x}_{bc}^j; \theta^{(i)})\right] \right)^2, \\
        \mathcal{L}_{r}(\theta) &= \frac{1}{N_{r}} \sum_{j=1}^{N_{r}} \left( \frac{\partial}{\partial t} \left[ \sum_{i=1}^{N} \omega_i(\mathbf{x}_{r}^j  )u_{\text{MHPKM}}^{(i)}(\mathbf{x}_{r}^j; \theta^{(i)}) \right] + \mathcal{N} \left[ \sum_{i=1}^{N} \omega_i(\mathbf{x}_{r}^j  )u_{\text{MHPKM}}^{(i)}(\mathbf{x}_{r}^j; \theta^{(i)}) \right] \right)^2.
    \end{aligned}
\end{equation*}%
Here, $\{\mathbf{x}_{ic}^j\}_{j=1}^{N_{ic}}$, $\{\mathbf{x}_{bc}^j\}_{j=1}^{N_{bc}}$ and $\{\mathbf{x}_{r}^j\}_{j=1}^{N_{r}}$ are the sample points for initial conditions (ICs), boundary conditions (BCs) and residual, respectively. $\lambda_{ic}$, $\lambda_{bc}$ and $\lambda_{r}$ in \eqref{eq:loss_function} are weight parameters to balance the contributions of each loss. The choice of these weights can significantly influence the training process and final performance of the network.

\begin{figure}[htbp]
    \centering 
        \includegraphics[scale=0.2]{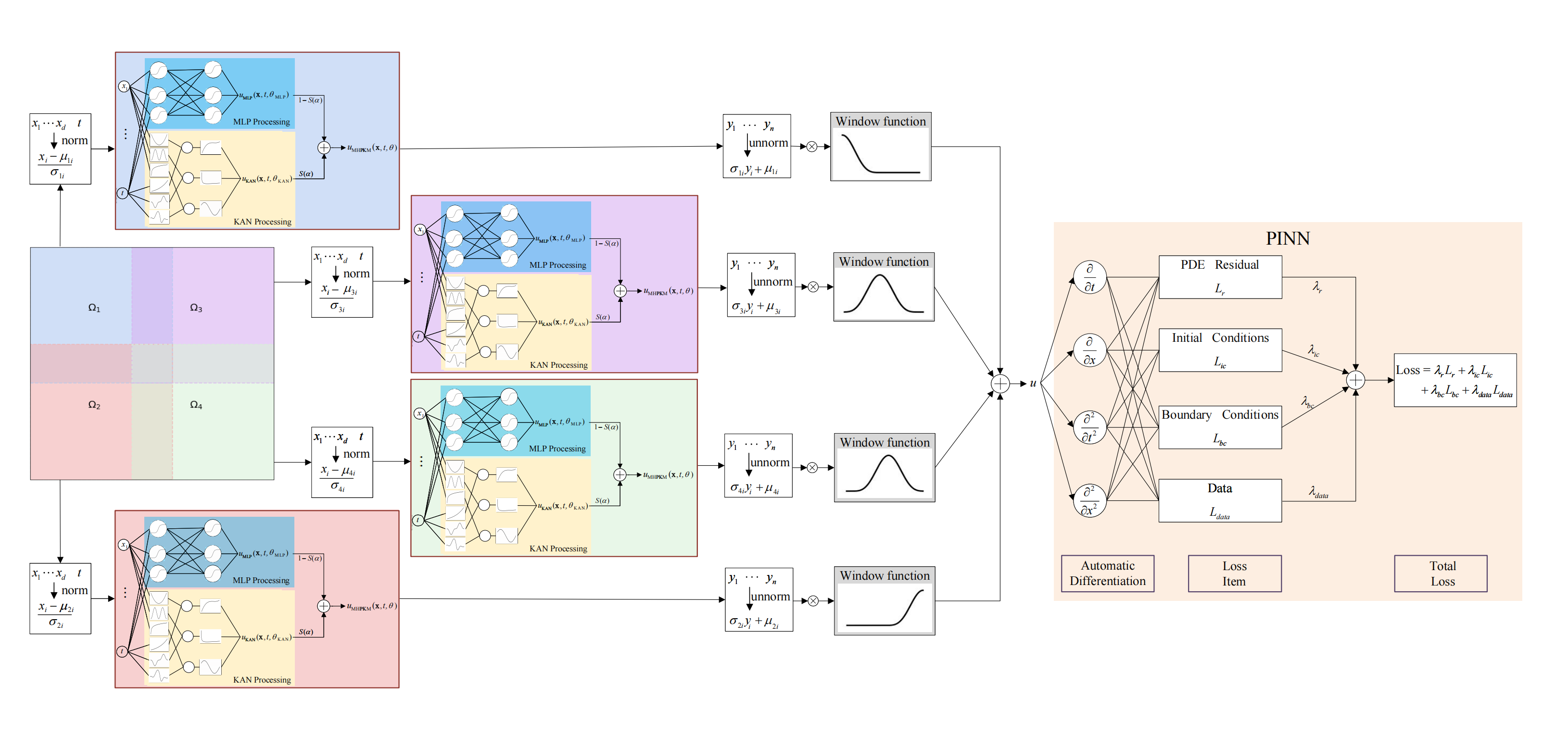}
   \caption{Graphical abstract of the modified HPKM-PINN with overlapping domain decomposition.}
    \label{Workflow Diagram}
\end{figure}

To resolve gradient competition among different loss terms in the loss function of the modified HPKM architecture $\eqref{eq:loss_function}$, we adopt a hard-constrained approach \cite{43}. By enforcing boundary and initial conditions through differential operator embedding, we can eliminate the corresponding loss terms and enable training with the simpler formulation

\begin{equation*}
\begin{aligned}    
    \mathcal{L}(\theta) = &\frac{1}{N_{r}} \sum_{j=1}^{N_{r}} \left( \frac{\partial}{\partial t} \left[ \sum_{i=1}^{N} \mathcal{C}[  \omega_i(\mathbf{x}_r^j)u_{\text{MHPKM}}^{(i)}(\mathbf{x}_{r}^j; \theta^{(i)})] \right] + \mathcal{N} \left[ \sum_{i=1}^{N}  \mathcal{C}[ \omega_i(\mathbf{x}_r^j)u_{\text{MHPKM}}^{(i)}(\mathbf{x}_{r}^j; \theta^{(i)})] \right] \right)^2,
\end{aligned}
\end{equation*}
where the operator $\mathcal{C}$ enforces hard constraints. The detailed computational procedure of our modified HPKM-PINN with overlapping domain decomposition is illustrated in Figure $\ref{Workflow Diagram}$ and Algorithm $\ref{Algorithm1}$.

\section{Numerical experiments}
In this section, we present some numerical experiments of the modified HPKM-PINN with overlapping domain decomposition for solving partial differential equations. Compared with a single-network baseline, the proposed approach demonstrates both superior efficiency and stronger generalization. We first assess the performance of various S-shaped functions in the network through tests on low-frequency and high-frequency problems. These tests demonstrate the modified HPKM architecture’s superior capability in handling multi-frequency problems. Subsequently, benchmark comparisons against standard MLP and KAN architectures further validate the proposed method's enhanced robustness and generalization performance.

In the following, all benchmark networks are trained with Adam at a learning rate of $1\times10^{-4}$. For comparison, we employ a 32-neuron MLP ($\tanh$ activation function) and a 32-neuron KAN (Fourier series activation function with $K=4$), while the modified HPKM-PINN uses only 16 neurons per layer for both its MLP and KAN components. To reduce optimization complexity, hard constraints are employed in all numerical examples to strictly enforce all boundary and initial conditions. The number of subdomains and the distribution of sampling points are determined through analysis according to problem complexity. For simplicity, we use a uniform overlap ratio $\delta$ for all regions. The performance of the networks is assessed by the normalized $L_{2}$ error, which is defined as 
\begin{equation}\label{eq:placeholder}
L_2(\tilde{u}; u) = \frac{\sqrt{\sum_{\mathbf{x}_i\in I} \left| u(\mathbf{x}_i) - \tilde{u}(\mathbf{x}_i;\theta) \right|^2}}{\sqrt{\sum_{\mathbf{x}_i\in I} \left| u(\mathbf{x}_i) \right|^2}},      
\end{equation}%
where $u$ denotes the reference solution, $\tilde{u}$ represents the prediction of the neural network, and $\mathbf{x}_i \in I = \{\mathbf{x}_1, \mathbf{x}_2, \dots\}$ is the sampling points.

\begin{algorithm}[H]
\caption{Modified HPKM-PINN with Overlapping Domain Decomposition}
\label{Algorithm1}
\begin{algorithmic}[1]
\State \textbf{Input:} Learning rate $\eta$, domain $\Omega$, number of subdomains $N$, sampling points $I$
\State Decompose $\Omega$ into $N$ overlapping subdomains $\{\Omega_i\}_{i=1}^N$
\State Initialize local modified HPKM networks $\{u^{(i)}(\mathbf{x};\theta^{(i)})\}_{i=1}^N$ parameterized by $\{\theta_{\text{KAN}}^{(i)},\theta_{\text{MLP}}^{(i)}\}_{i=1}^N$
\State Initialize trainable weighted parameters $\{\alpha^{(i)}\}_{i=1}^N = 0$ 
\State \# The set of all trainable parameters in the modified HPKM-PINN is $\theta = \{\theta_{\text{KAN}}^{(i)},\theta_{\text{MLP}}^{(i)}, \alpha^{(i)}\}_{i=1}^N $
\While{not converged}    
    \For{$i = 1$ to $N$} \Comment{Process each subdomain in parallel}
        \For{each point $\mathbf{x} \in \Omega_i$}
            \State $\tilde{\mathbf{x}}^{(i)} \gets \text{Normalize}(\mathbf{x})$ \Comment{Map to $[0,1]^d$}
            \State $\tilde{u}^{(i)} \gets u^{(i)}(\tilde{\mathbf{x}}^{(i)};\theta^{(i)})$ \Comment{Local network forward pass}
            \State $u^{(i)} \gets \text{Unnormalize}(\tilde{u}^{(i)})$ \Comment{Restore to physical scale}
        \EndFor
    \EndFor

    \State Hard constraints: $\hat{u}^{(i)} = \mathcal{C} u^{(i)} $
    
    \State Construct global solution: $\tilde{u}(\mathbf{x},\theta) = \sum_{i=1}^N \omega_i(\mathbf{x}) \cdot \hat{u}^{(i)}$
    
    \State Total loss: $\mathcal{L}(\theta) = \lambda_r\mathcal{L}_r$

\State Update the parameters $\theta $ via gradient descent

     $$\theta_{n+1} = \theta_{n} - \eta \nabla_{\theta_{n}} \mathcal{L}(\theta_n)$$

\EndWhile

\State \Return  $\tilde{u}(\mathbf{x},\theta)$
\end{algorithmic}
\end{algorithm}

\subsection{Helmholtz equation}
The Helmholtz equation is an elliptic partial differential equation that plays a role in the modeling of wave propagation phenomena \cite{zengjiaHelmhotz}. This governing equation has extensive applications in seismology, electromagnetic radiation, acoustics, and various physical systems. And it is equivalent to the wave equation under the assumption of a single frequency. 
In this example, we consider the following Helmholtz equation with
homogeneous Dirichlet boundary condition 
\begin{equation}   \label{Hemholtz equation}
\begin{aligned}
    &\Delta u(x,y)+\kappa^{2} u(x,y)=h(x,y),~~~&&(x,y)\in \Omega= (-1,1)^2,\\
&u(x,y)=0, &&(x,y) \in \partial \Omega,
\end{aligned}
\end{equation}
where $u(x, y)$ denotes a spatially distributed field quantity, $\kappa$ is the wave number. $\Delta$ is the Laplacian operator, $h(x, y)$ is a forcing term that introduces inhomogeneity into the system. Specifically, $h$ is defined as
\begin{equation*}
 h(x,y) = (1-2\omega^2\pi^2)\sin{(\omega \pi x)}\sin{(\omega \pi y)},
\end{equation*}
where $\omega$ is a constant.
The exact solution of the problem ($\ref{Hemholtz equation}$) is given by 
\begin{equation*}
u(x,y)=\sin(\omega \pi x)\sin(\omega \pi y).
\end{equation*}
The scale discrepancy between the residual loss and boundary loss leads to training instability. To solve this problem and ensure that the network satisfies the specified boundary conditions, we enforce hard constraints in the following form 
\begin{equation*}
\tilde{u}(x, y) = \tanh(\frac{x+1}{\sigma})\tanh(\frac{1-x}{\sigma})\tanh(\frac{y+1}{\sigma})\tanh(\frac{1-y}{\sigma})NN(x,y;\theta), 
\end{equation*}
where $NN(x, y; \theta)$ is the output of the neural network, $\sigma$ is a scaling factor that controls the steepness of the tanh functions. In this example, we set $\sigma = 0.2 $ and $\kappa = 1.0$. The following experiments are conducted under uniform parameter settings: learning rate $1\times10^{-4}$, domain decomposition of $20\times 20$ subdomains, and $260\times 260$ sampling points. Based on our tests, the overlap ratio $\delta$ is selected as $3.3$. This value achieves a good balance between computational efficiency and accuracy.

\subsubsection{The test of different S-shaped functions}
\label{subsec:sigmoidal}

To investigate the impact of the S-shaped weighting function selection on the output of the modified HPKM-PINN with trainable parameter $\alpha$, we conduct comparative tests of multiple S-shaped functions $S(\alpha)$, including that
\begin{equation}
\label{sigmoidal_function}
\begin{aligned}
&\text{Sigmoid function:}~~S(\alpha) = \dfrac{1}{1 + e^{-\alpha}}, 
~&&       
\text{Tanh scaled function:}~~S(\alpha) = \dfrac{\tanh(\alpha) + 1}{2}, 
\\
&\text{Arctangent scaled function:}~~S(\alpha) = \dfrac{\arctan(\alpha)}{\pi} + \frac{1}{2},
~&&
\text{Softsign scaled function:}~~S(\alpha) = \frac{1}{2} \left( \dfrac{\alpha}{1 + |\alpha|} + 1 \right),
\\        
&\text{Algebraic sigmoid function:}~~S(\alpha) = \frac{1}{2} \left( \dfrac{\alpha}{\sqrt{1 + \alpha^2}} + 1 \right),
~&&        
\text{Clip function:}~~S(\alpha) = \begin{cases} 
            1 & \text{if } \alpha > 1, \\ 
            0 & \text{if } \alpha < 0, \\ 
            \alpha & \text{otherwise}.
        \end{cases}
\end{aligned}
\end{equation}

We remark that when employing the clip function, the parameter $\alpha$ is trained directly without utilizing a S-shaped function $S(\alpha)$ in its learning dynamics. In addition to the considered S-shaped functions in \eqref{sigmoidal_function}, we also examine the unweighted HPKM approach that the output is obtained by summing the outputs of the MLP and KAN directly.

First, we consider the case when $\omega=1$, i.e. the exact solution to the problem (\ref{Hemholtz equation}) is 
\begin{equation*}
u(x,y)=\sin( \pi x)\sin( \pi y).
\end{equation*}

\begin{figure}[h]
    \centering
    \begin{minipage}[b]{0.3 \textwidth}  
        \centering
        \includegraphics[width=1.21\textwidth]{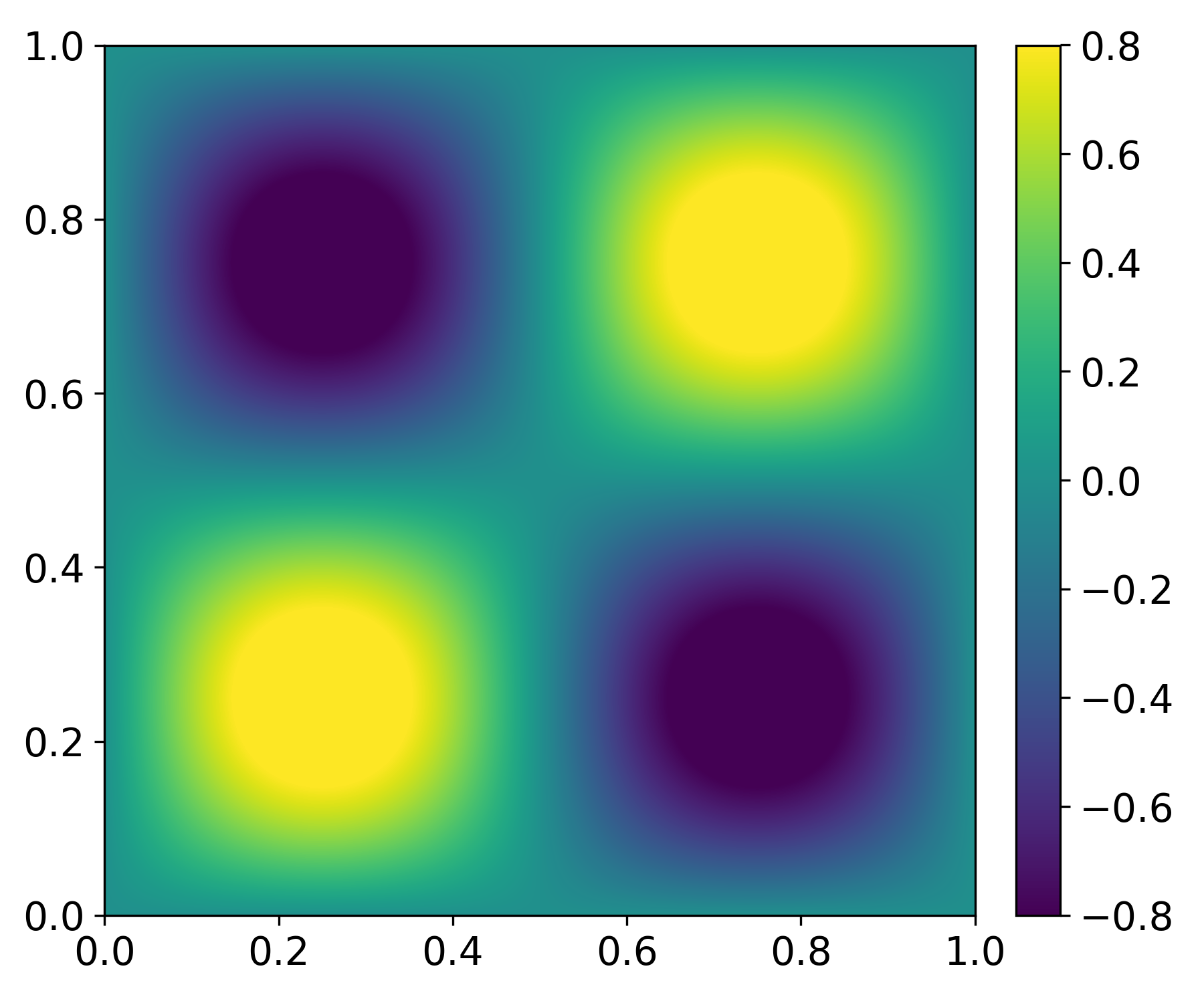}   
        \subcaption{}
    \end{minipage}\hspace{0.105 \textwidth} 
    \begin{minipage}[b]{0.3 \textwidth}  
        \centering
        \includegraphics[width=1.21\textwidth]{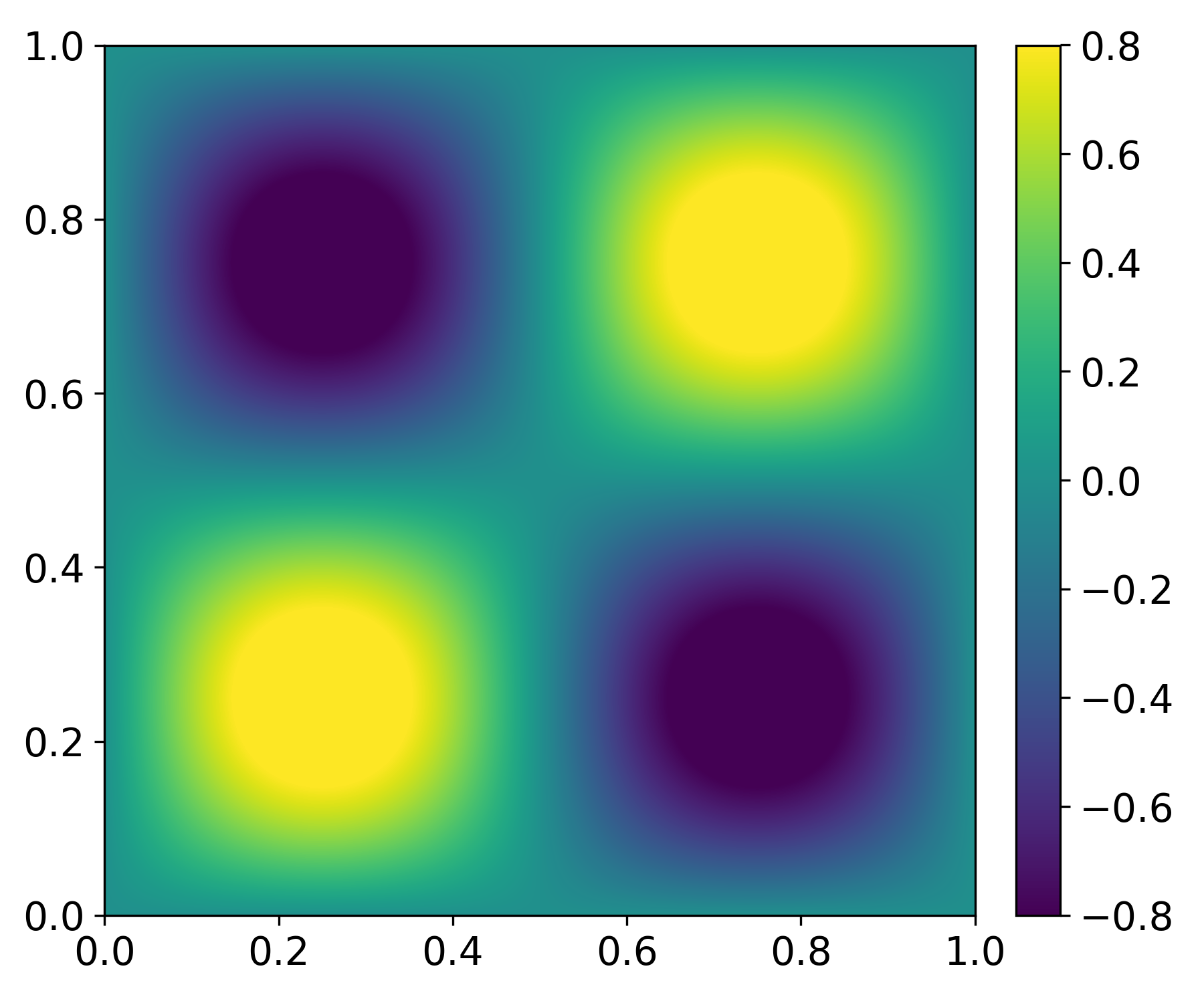}   
        \subcaption{}
    \end{minipage}
    
    \vspace{0.2cm} 
    \begin{minipage}[b]{0.3\textwidth}
        \centering
        \includegraphics[width=1.1 \textwidth]{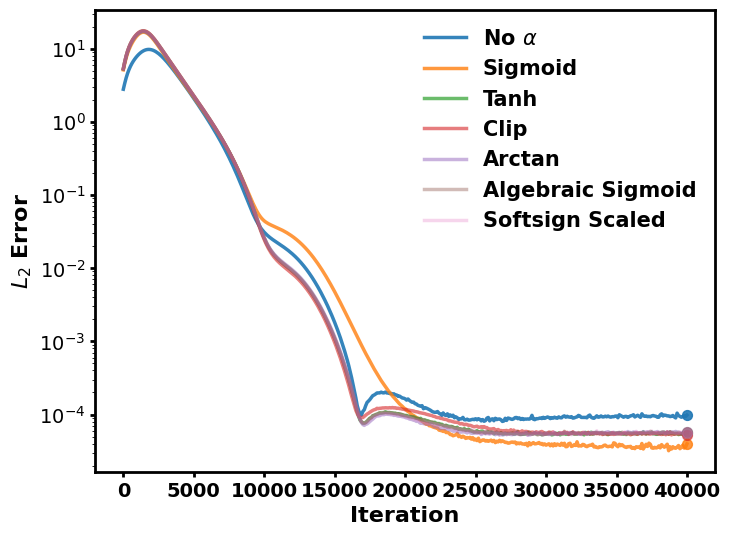}   
        \subcaption{}
    \end{minipage}\hspace{0.08\textwidth}
    \begin{minipage}[b]{0.3\textwidth}
        \centering
        \includegraphics[width=1.1\textwidth]{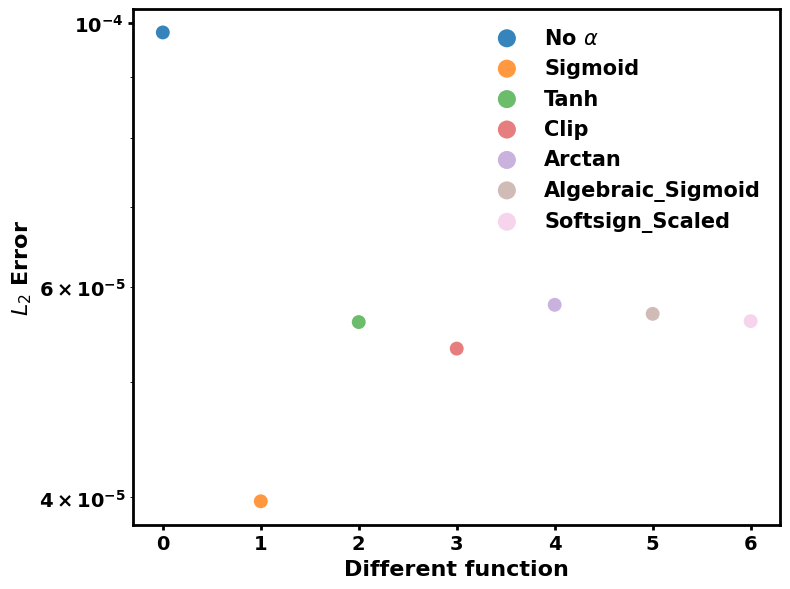}   
        \subcaption{}
    \end{minipage} 
    
    \caption{Helmholtz equation ($\omega=1$): (a) Exact solution; (b) Predicted solution; (c) The convergence curve of $L_2$ error with different S-shaped functions; (d) The final $L_2$ error after completing the training with different S-shaped functions.}  
    \label{darcylaw_l23}
\end{figure}

As shown in Figure \ref{darcylaw_l23}, the exact solution exhibits low-frequency characteristics in (a), and the predicted solution with a sigmoid $S(\alpha)$ is presented in (b). The dynamic evolution of the $L_2$ error during training for different S-shaped functions described in (3.3) is illustrated in (c), and (d) compares the final $L_2$ errors after completing training. The results indicate that all networks are capable of solving this problem well. When the S-shaped function is selected as the sigmoid function, the final $L_2$ error reaches approximately $4\times10^{-5}$, while for other functions in (\ref{sigmoidal_function}), the error remains around $5.6\times10^{-5}$. Moreover, compared to the unweighted MLP-KAN direct summation approach whose error reaches $9.8\times10^{-5}$ (see the `No $\alpha$' curve in Figure \ref{darcylaw_l23} (c)-(d), all weighted networks show improved accuracy, which confirms the role of $S(\alpha)$ in enhancing predictive performance.

\begin{figure}[h]
    \centering
    \begin{minipage}[b]{0.31\textwidth}  
    \includegraphics[width=1.1\textwidth]{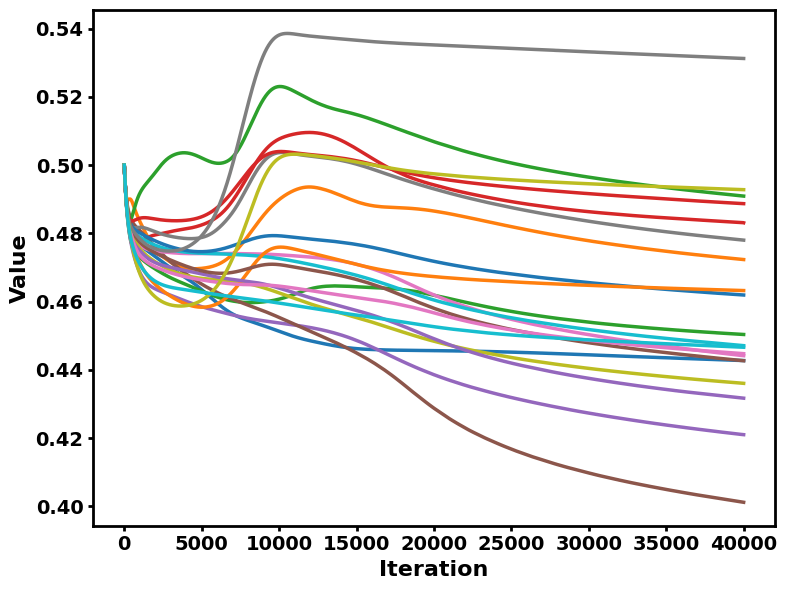}
        \centering
        \subcaption{$S(\alpha)$: Clip function}
    \end{minipage}
    \hfill
    \begin{minipage}[b]{0.31\textwidth}
        \includegraphics[width=1.1\textwidth]{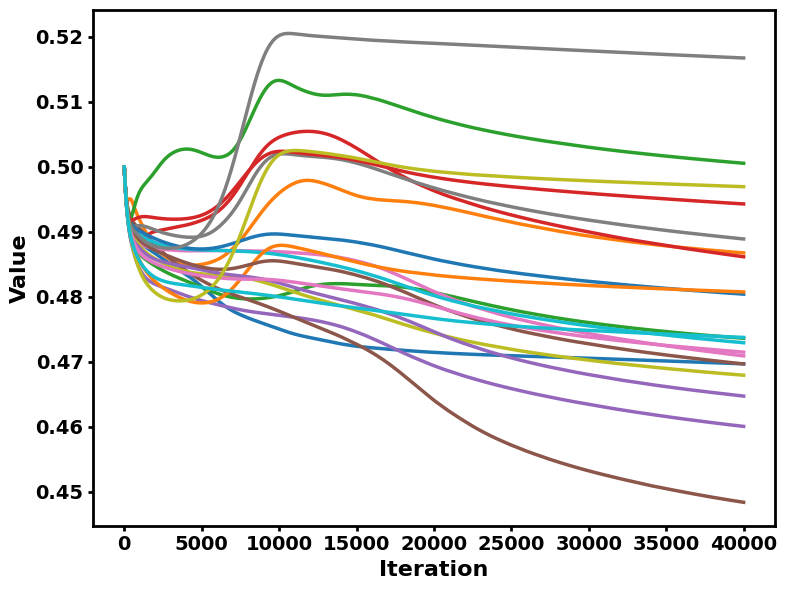}   
        \centering
       \subcaption{$S(\alpha)$: Tanh scaled function}
    \end{minipage}
    \hfill
    \begin{minipage}[b]{0.31\textwidth}
        \includegraphics[width=1.1\textwidth]{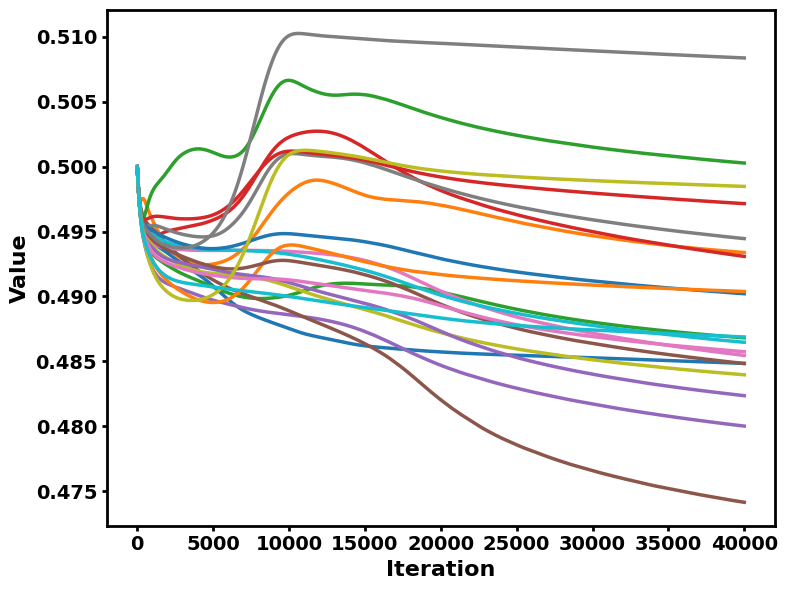}    
        \centering
        \subcaption{$S(\alpha)$: Sigmoid function}
    \end{minipage}
   \caption{ 20 random evolution curves of the weights $S(\alpha)$ for Helmholtz equation ($\omega=1$).} 
    \label{fig:darcy_results}
\end{figure}

To investigate the evolution of parameter $S(\alpha)$, we visualize its dynamic evolution throughout the optimization process. Figure \ref{fig:darcy_results} shows the evolutionary patterns of $S(\alpha)$ in 20 randomly selected subdomains for three functions: the clip function, the tanh scaled function, and the sigmoid function. The evolution of $S(\alpha)$ exhibits a nearly consistent trend across the three functions. Furthermore, the introduction of S-shaped functions helps constrain the variation range of the convex combination coefficients in the modified HPKM architecture, thereby aiding the network in identifying more effective coefficients. In particular, the value of $S(\alpha)$ shows a declining trend, reflecting the increased dependence of the model on the output of the MLP branch. This behavior validates the theoretical expectation that MLP is effective to capture global features for the low-frequency problems considered in this test.    

We further explore the influence of the S-shaped function on the high-frequency oscillatory characteristics of the solutions to the Helmholtz equation. We consider $\omega = 32$ in  ($\ref{Hemholtz equation}$), where the exact solution is
\begin{equation*}
u(x,y)=\sin(32 \pi x)\sin(32 \pi y).
\end{equation*}

\begin{figure}[h]
    \centering
    \begin{minipage}[b]{0.29\textwidth}
        \centering
        \includegraphics[width=1.21\textwidth]{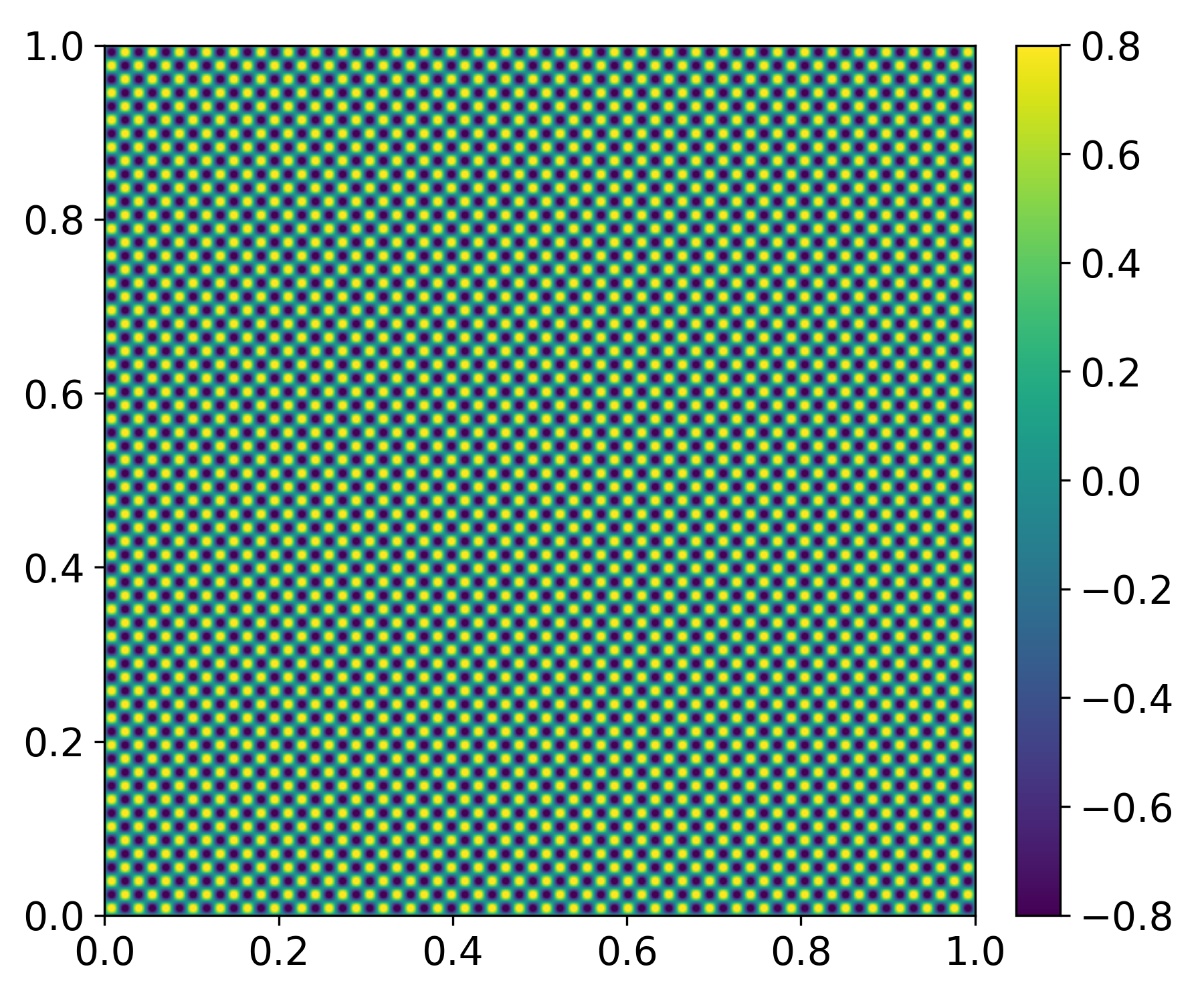}   
         \subcaption{}
    \end{minipage}\hspace{0.11\textwidth}
    \begin{minipage}[b]{0.29\textwidth}
        \centering
        \includegraphics[width=1.21\textwidth]{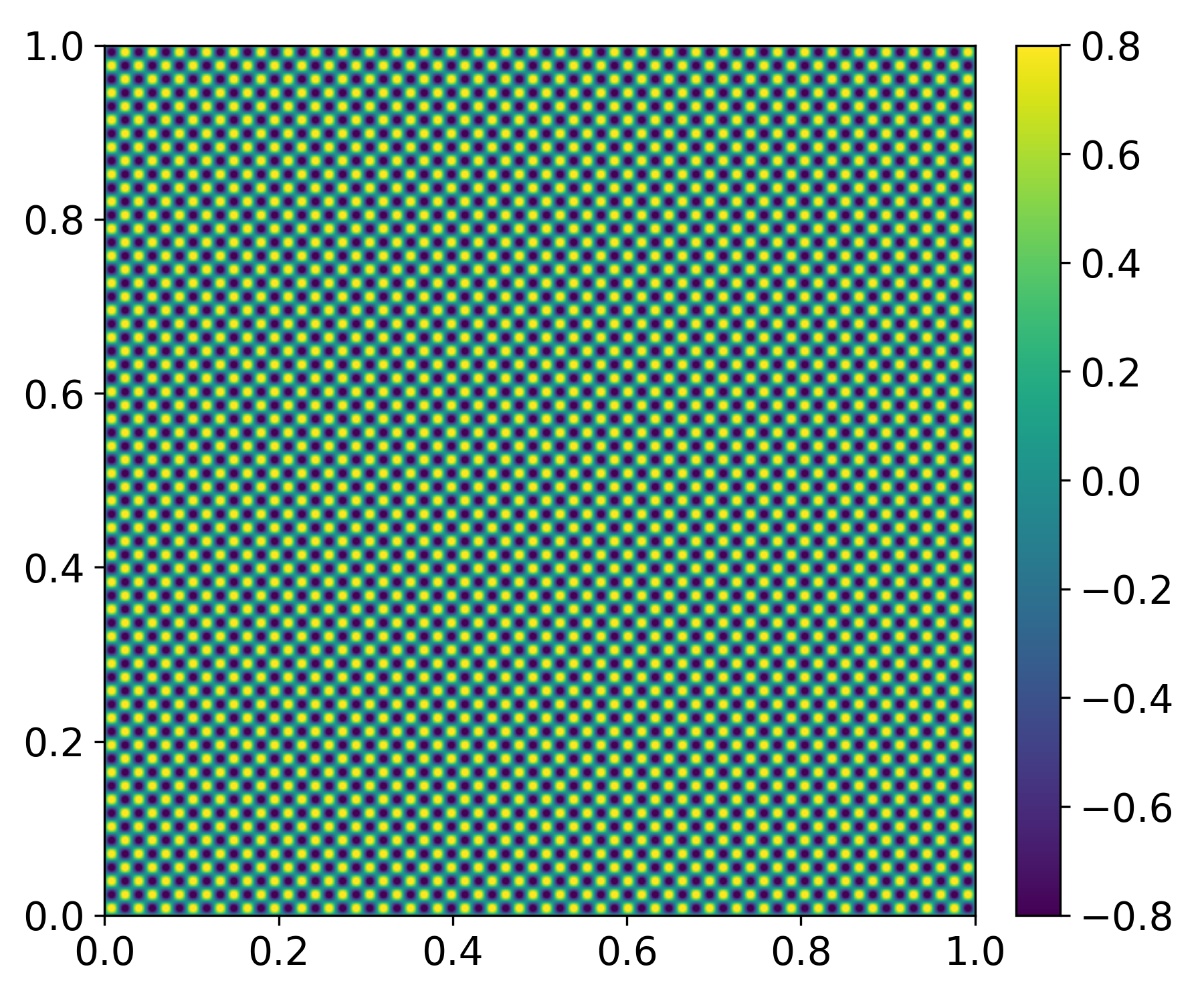}  
        \subcaption{}
    \end{minipage}
    \vspace{0.1cm} 
    \begin{minipage}[b]{0.3\textwidth}
        \centering
        \includegraphics[width=1.1\textwidth]{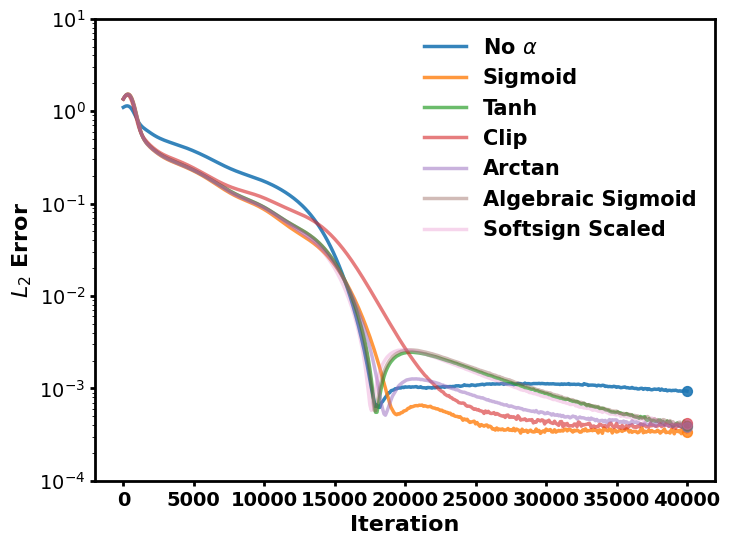}   
       \subcaption{}
    \end{minipage}\hspace{0.08\textwidth}
    \begin{minipage}[b]{0.3\textwidth}
        \centering
        \includegraphics[width=1.1\textwidth]{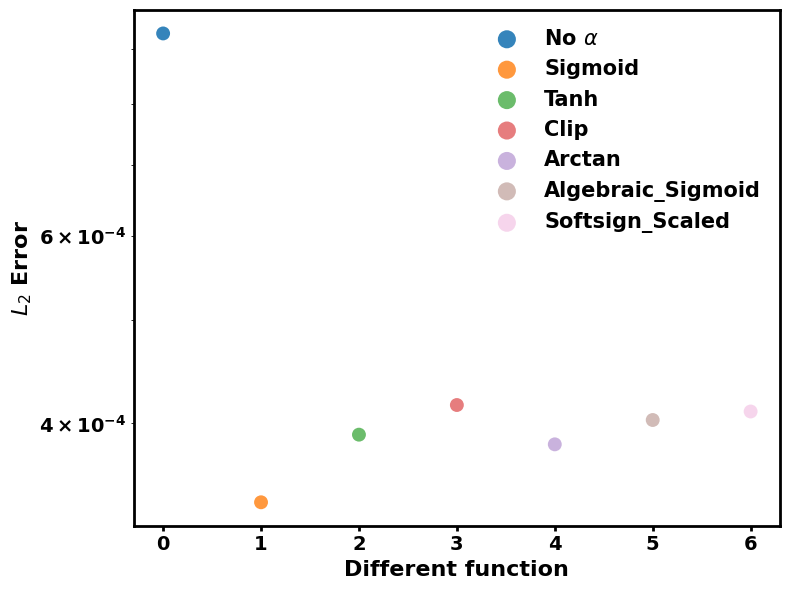}   
      \subcaption{}
    \end{minipage}
    
    \caption{Helmholtz equation ($\omega=32$): (a) Exact solution; (b) Predicted solution; (c) The convergence curve of $L_2$ error with different S-shaped functions; (d) The final $L_2$ error after completing the training with different S-shaped functions.}  
    \label{Helmholtz_equation_l23}
\end{figure}

The results are illustrated in Figure $\ref{Helmholtz_equation_l23}$: the exact and predicted solutions (using a sigmoid $S(\alpha)$) are shown in (a) and (b), respectively; the $L_2$ error results across different training stages when applying various S-shaped functions to $\alpha$ is presented in (c); and (d) summarizes the final $L_2$ error comparison after completing training. These results demonstrate that for high-frequency problems ($\ref{Hemholtz equation}$), the modified HPKM-PINN where the S-shaped function is considered, achieves a higher convergence accuracy compared to directly training the parameter $\alpha$ using the clip function. Furthermore, the modified HPKM architecture incorporating the adaptive weighting functions $S(\alpha)$ in (\ref{sigmoidal_function}) show superior accuracy over the unweighted MLP-KAN approach (see the `No $\alpha$' curve in Figure \ref{Helmholtz_equation_l23} (c)-(d)). These findings indicate that the adaptive weighting function $S(\alpha)$ enhances both model accuracy and training stability, and is beneficial in solving high-frequency problems.

\begin{figure}[htbp]
    \centering
    \begin{minipage}[b]{0.31\textwidth}
        \includegraphics[width=1.1\textwidth]{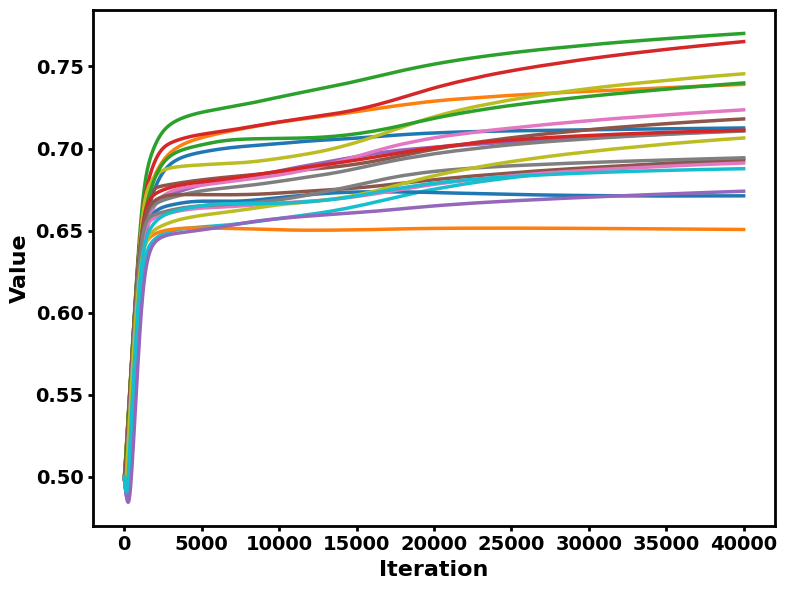}
        \centering
        \subcaption{$S(\alpha)$: Clip function}
    \end{minipage}
    \hfill
    \begin{minipage}[b]{0.31\textwidth}
        \includegraphics[width=1.1\textwidth]{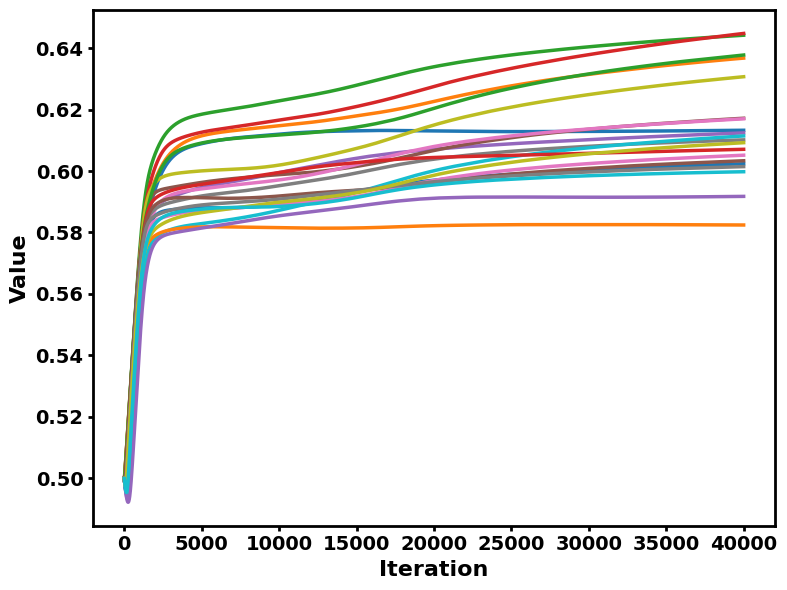}
        \centering
        \subcaption{$S(\alpha)$: Tanh scaled function}
    \end{minipage}
    \hfill
    \begin{minipage}[b]{0.31\textwidth}
        \includegraphics[width=1.1\textwidth]{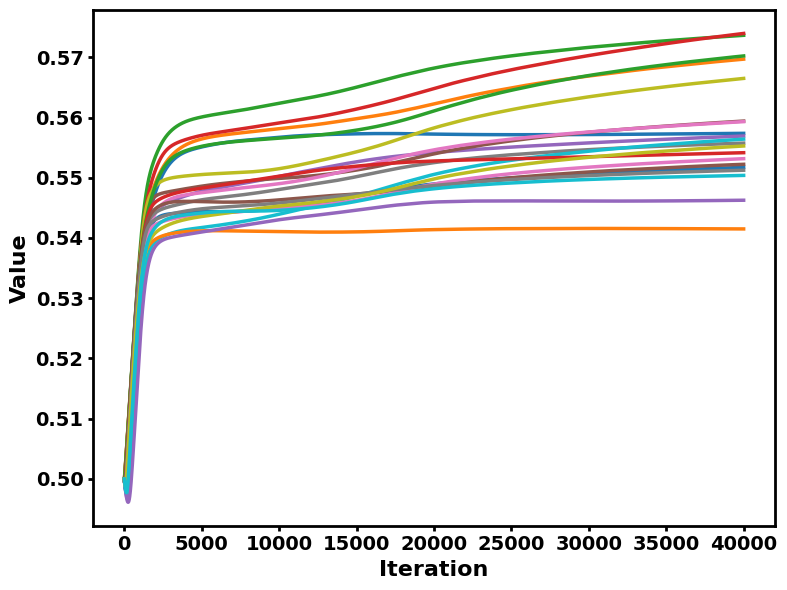}
        \centering
        \subcaption{$S(\alpha)$: Sigmoid function}
    \end{minipage}
    \caption{ 20 random evolution curves of the weights $S(\alpha)$ for Helmholtz equation ($\omega=32$).} 
    \label{fig:helmholtz_results}
\end{figure}

Similarly, we analyze the evolution of $S(\alpha)$ in high-frequency Helmholtz problems. Figure $\ref{fig:helmholtz_results}$ shows the evolutionary patterns of $S(\alpha)$ in 20 randomly selected subdomains. It should be noted that the unpresented subdomains exhibit consistent trends with the selected representatives. Figure $\ref{fig:helmholtz_results}$ (a) illustrates the distribution of $S$ under the clip function $S(\alpha)$ in the chosen subdomains, Figure $\ref{fig:helmholtz_results}$ (b) shows the evolutionary process of the tanh scaled function $S(\alpha)$, and Figure $\ref{fig:helmholtz_results}$ (c) documents the dynamic variation of the sigmoid function $S(\alpha)$ during training. The evolution of the value $S(\alpha)$ reveals a consistent pattern across all three figures. The introduction of S-shaped functions reduces the variation range of the weight parameter $S(\alpha)$, indicating that only minor weight adjustments are needed to achieve acceptable accuracy. Furthermore, the upward trend observed in $S(\alpha)$ across all experiments demonstrates that the model favors the output from the KAN branch during the learning process. This phenomenon aligns with our theoretical expectations: the KAN network's enhanced local feature extraction capability leads to superior performance in solving high-frequency problems.

The phenomena observed in Figure $\ref{fig:darcy_results}$ and Figure $\ref{fig:helmholtz_results}$ align with the intrinsic properties of the model components. The MLP architecture maintains stable performance across the full frequency spectrum through its broadband learning capability, excelling in global feature extraction, and the KAN architecture demonstrates superior nonlinear feature extraction, particularly effective in capturing localized detailed features. In addition, as shown in Figure $\ref{darcylaw_l23}$ (d) and $\ref{Helmholtz_equation_l23}$ (d), employing the sigmoid function to train the weighting function $S(\alpha)$ not only maintains the weight range of 0 to 1, but also improves the network's accuracy. We select the sigmoid function as the S-shaped function to optimize the parameter $\alpha$ due to its superior performance in our tests.

\subsubsection{Modified HPKM vs MLP and KAN}

\begin{figure}[htbp]
\centering
\begin{minipage}[b]{0.3\textwidth}
  \centering
  \includegraphics[width=0.9\linewidth]{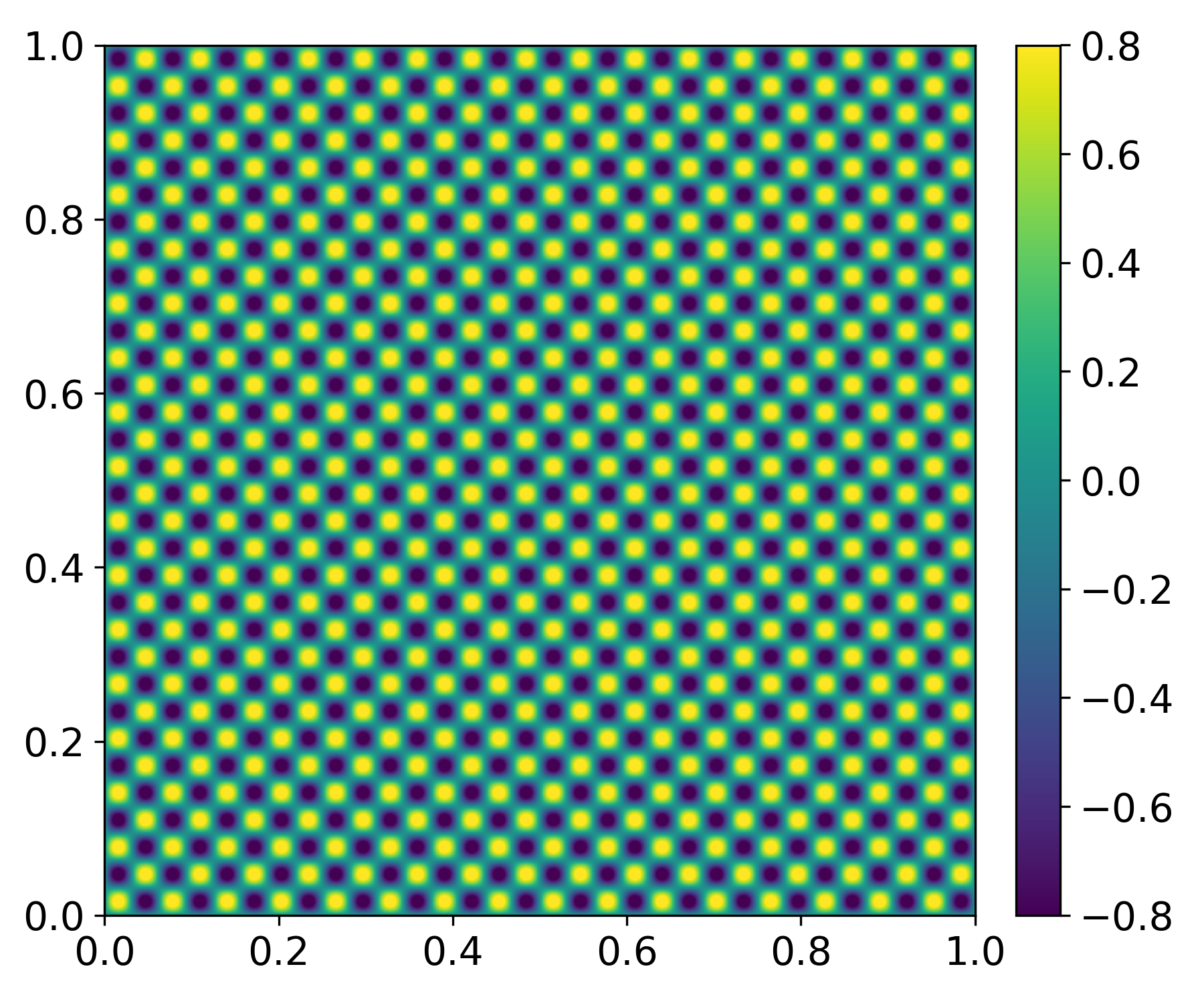}
  \vspace{0.2pt}
  \footnotesize Exact solution $(\omega=16)$
\end{minipage}
\hfill
\begin{minipage}[b]{0.3\textwidth}
  \centering
  \includegraphics[width=0.9\linewidth]{Pic/function_omega_32.png}
  \vspace{0.2pt}
  \footnotesize Exact solution $(\omega=32)$
\end{minipage}
\hfill
\begin{minipage}[b]{0.3\textwidth}
  \centering
  \includegraphics[width=0.9\linewidth]{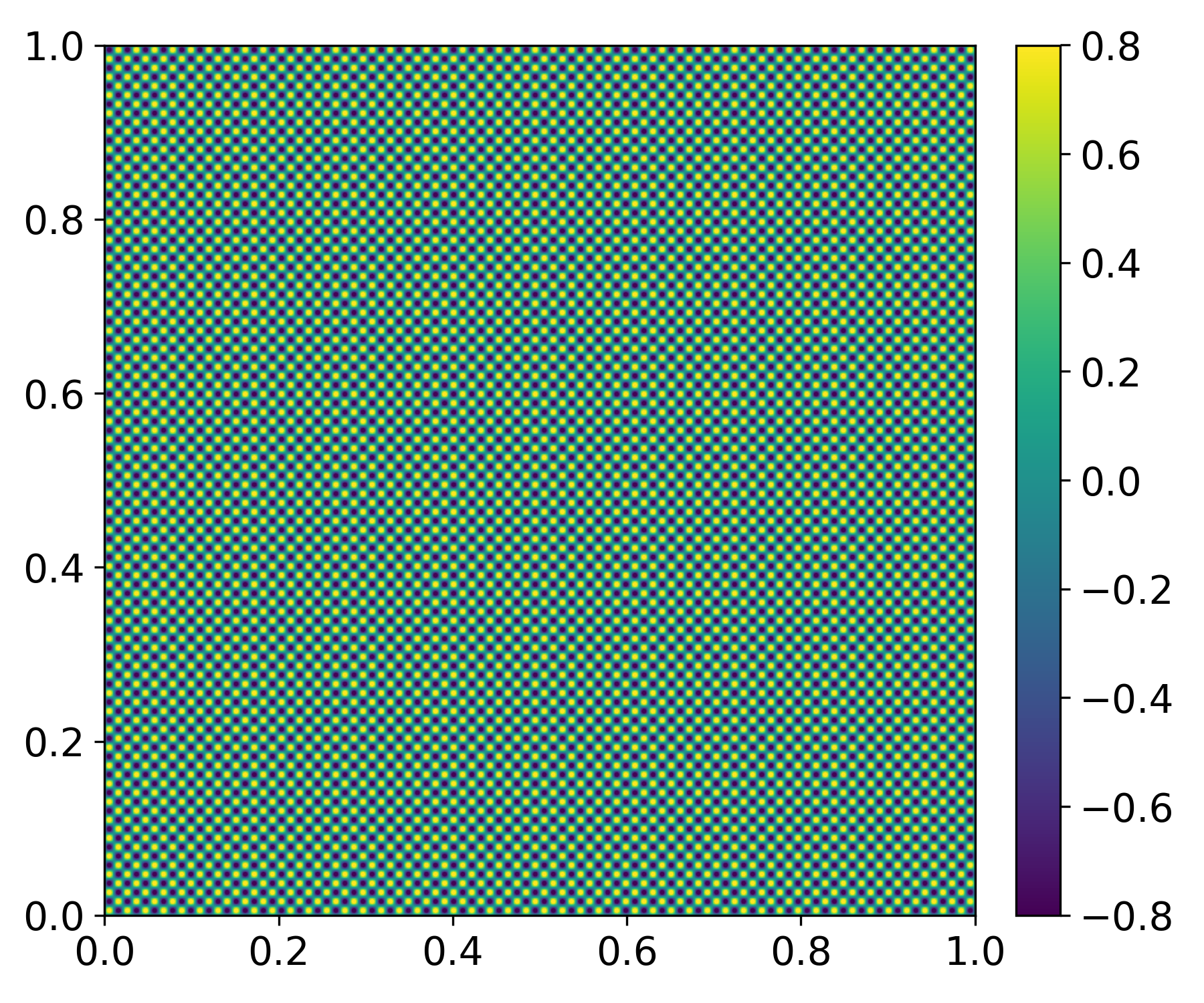}
  \vspace{0.2pt}
  \footnotesize Exact solution $(\omega=48)$
\end{minipage}

\vspace{2ex}  

\begin{minipage}[b]{0.3\textwidth}
  \centering
  \includegraphics[width=0.9\linewidth]{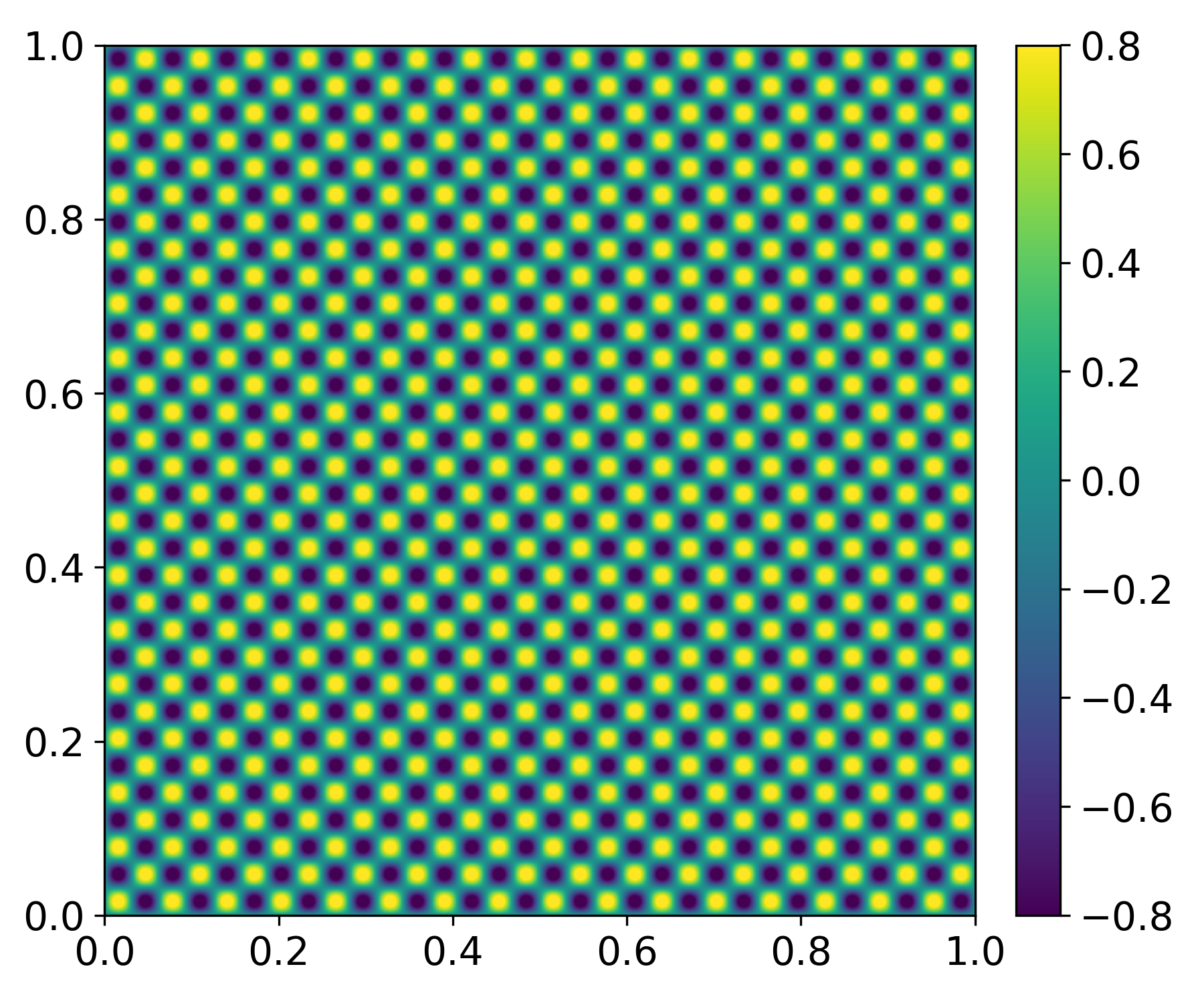}
  \vspace{1pt}
  \footnotesize Predicted solution $(\omega=16)$
\end{minipage}
\hfill
\begin{minipage}[b]{0.3\textwidth}
  \centering
  \includegraphics[width=0.9\linewidth]{Pic/predict_omega_32.png}
  \vspace{1pt}
  \footnotesize Predicted solution $(\omega=32)$
\end{minipage}
\hfill
\begin{minipage}[b]{0.3\textwidth}
  \centering
  \includegraphics[width=0.9\linewidth]{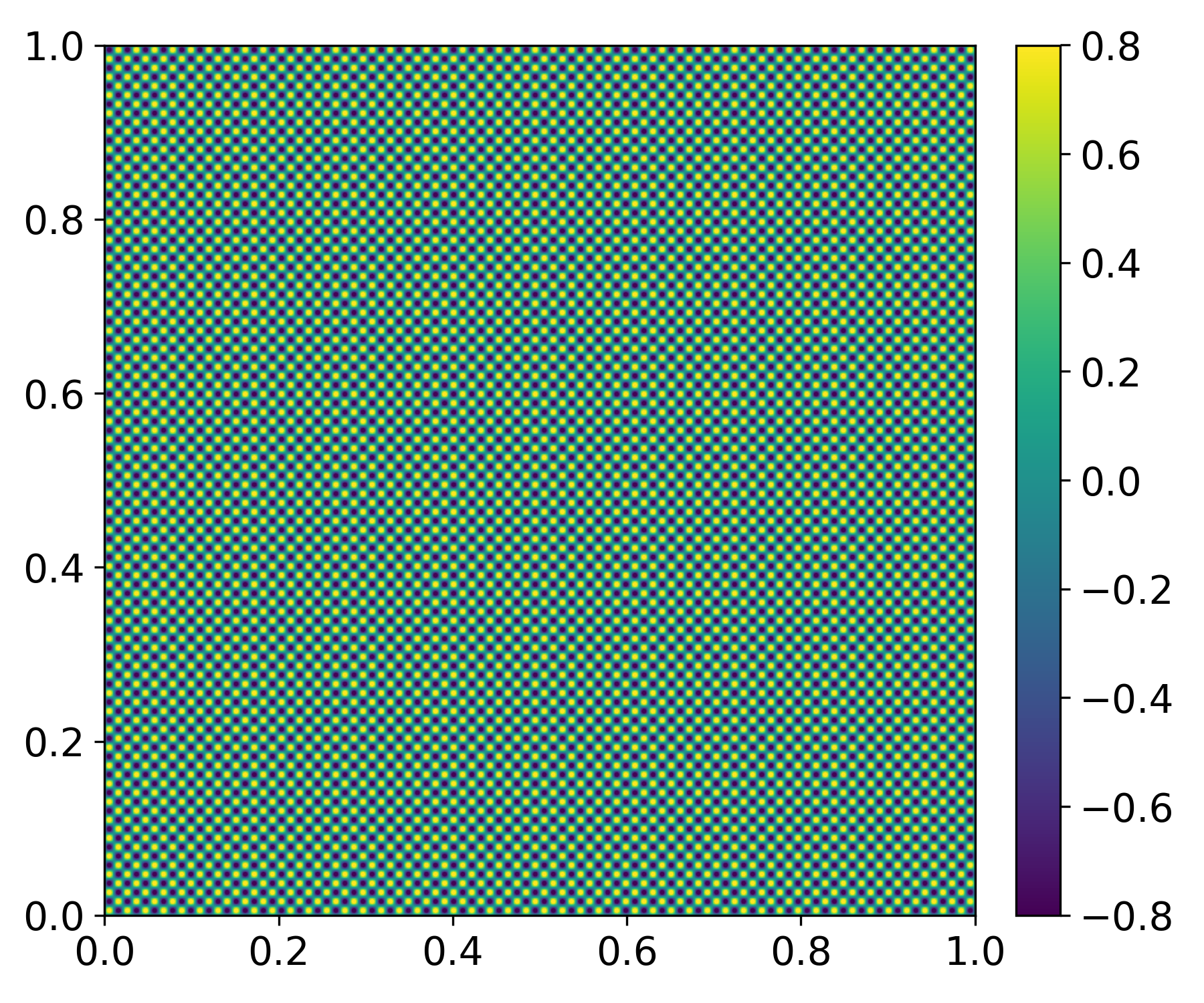}
  \vspace{1pt}
  \footnotesize Predicted solution $(\omega=48)$
\end{minipage}

\vspace{2ex}

\begin{minipage}[b]{0.3\textwidth}
  \centering
  \includegraphics[width=0.9\linewidth]{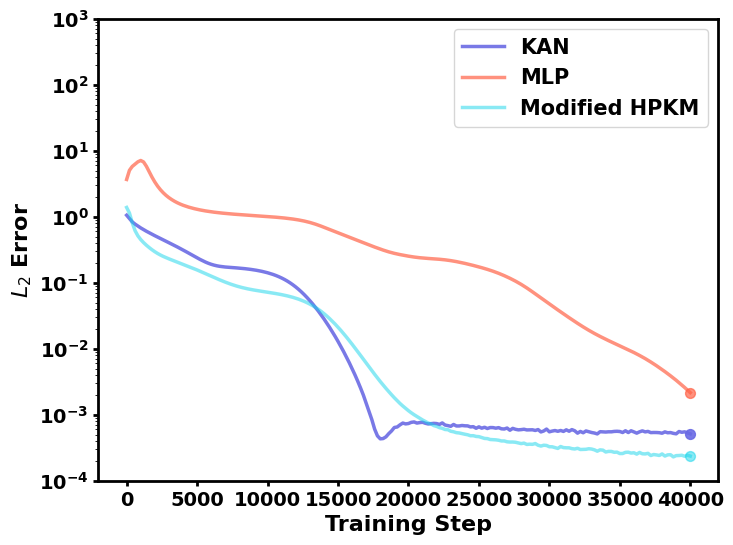}
  \vspace{1pt}
  \footnotesize $L_2$ error $(\omega = 16)$
\end{minipage}
\hfill
\begin{minipage}[b]{0.3\textwidth}
  \centering
  \includegraphics[width=0.9\linewidth]{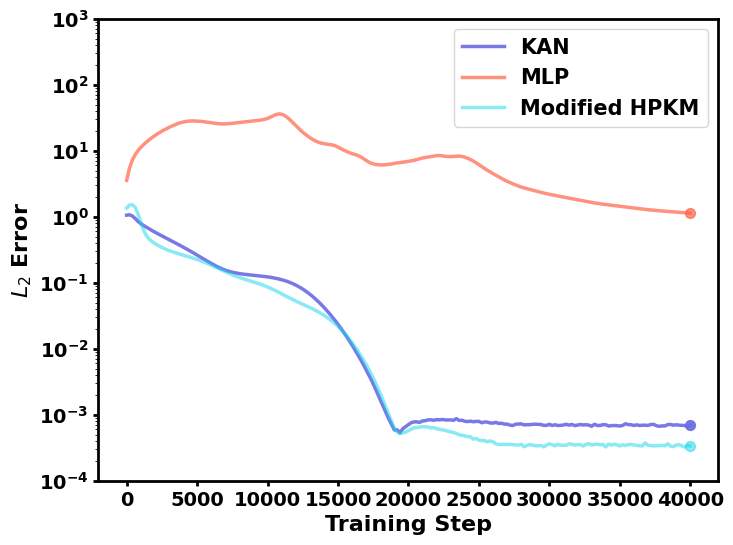}
  \vspace{1pt}
  \footnotesize $L_2$ error $(\omega = 32)$
\end{minipage}
\hfill
\begin{minipage}[b]{0.3\textwidth}
  \centering
  \includegraphics[width=0.9\linewidth]{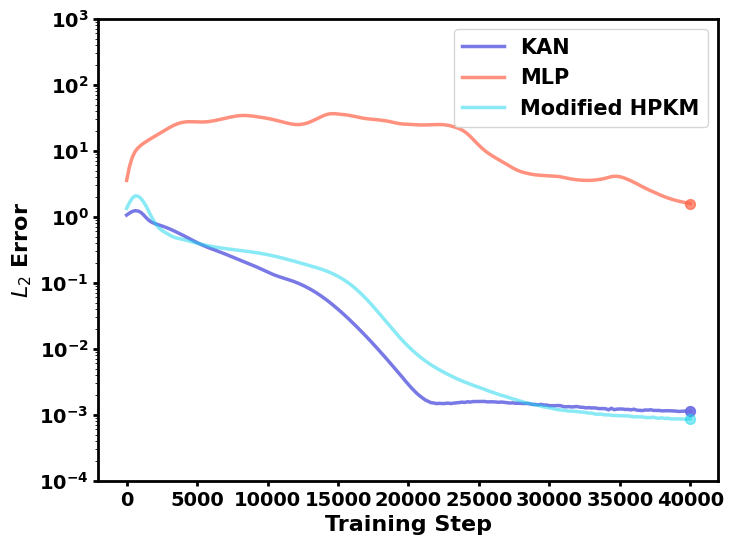}
  \vspace{1pt}
  \footnotesize $L_2$ error $(\omega = 48)$
\end{minipage}

\caption{Helmholtz equation with $\omega=16,32,48$. Top: Exact solutions. Middle: Predicted solutions obtained by the modified HPKM-PINN. Bottom: $L_2$ error comparison between KAN, MLP and the modified HPKM architectures.}
\label{results-Helmholtz-w-equation}
\end{figure}

In the following, employing a domain decomposition technique, we evaluate the performance of MLP-based PINN model, KAN-baed PINN model, and the modified HPKM-PINN model in solving the Helmholtz equation with varying frequency characteristics. The top and second lines of Figure $\ref{results-Helmholtz-w-equation}$ show that the exact solution and the predicted solution develop increasingly significant high-frequency oscillations as the parameter \(\omega\) grows ($\omega=16,32,48$). The result of $L_2$ error is shown that under the current hyperparameter configuration, the PINN based on MLP with domain decomposition technique shows excessive sensitivity to hyperparameter settings, requiring laborious fine-tuning to accommodate problem complexity. However, the modified HPKM architecture consistently demonstrates superior performance. This advantage stems from the the modified HPKM-PINN's KAN architecture and the adaptive dual network balancing mechanism through trainable parameter $S(\alpha)$. Furthermore, when the solution becomes increasingly oscillatory and complex as shown in the bottom line of Figure $\ref{results-Helmholtz-w-equation}$, the PINN based on MLP with domain decomposition technique demonstrates inferior performance compared to the other two models, but the KAN architecture with domain decomposition technique exhibits superior advantages for high-frequency problems. More importantly, by leveraging an adaptive weighting function $S(\alpha)$, the modified HPKM-PINN integrates the strengths of KAN and MLP, resulting in improved performance.

The modified HPKM architecture, which incorporates a convex combination of the outputs of KAN and MLP, not only improves model performance but also exhibits superior computational efficiency and lower training costs compared to the KAN architecture. Analysis of the convex combination weighting function $S(\alpha)$ reveals that in this case the modified HPKM-PINN model is more dependent on the KAN output during training. All the results in Figure $\ref{results-Helmholtz-w-equation}$ show that the modified HPKM architecture adapts the output contributions of the two architectures through the S-shaped function, integrating the strengths of MLP and KAN architectures to achieve robust performance across diverse complex problems. The $L_2$ errors presented in Table \ref{err_total_result} demonstrate that the modified HPKM-PINN model combined with domain decomposition addresses high-frequency problems more effectively and achieves higher accuracy. This result highlights the promise of the modified HPKM-PINN for addressing such challenges.

\subsection{Poisson equation}
The Poisson equation describes the spatial distribution of potential fields generated by given source terms, effectively capturing the interplay between diffusion processes and external forcing mechanisms. As a cornerstone of nonlinear mathematical physics, the Poisson equation finds widespread applications across engineering and applied mathematics. To investigate the performance of the modified HPKM-PINN with overlapping domain decomposition in handling multiscale phenomena and high-frequency features, we conduct numerical experiments on two-dimensional and high-dimensional Poisson equations with deliberately increased complexity.

Specifically, we consider the following form of the Poisson equation 
\begin{equation}  \label{d Poisson equation}
\begin{aligned}
-\Delta u(\mathbf{x}) &= h(\mathbf{x}), && \mathbf{x} \in \Omega = (0,1)^d, \\
u(\mathbf{x}) &= 0, && \mathbf{x} \in \partial \Omega.
\end{aligned}
\end{equation}

\subsubsection{Two-dimensional Poisson equation}
We first consider the two-dimensional case ($d=2$) of the Poisson problem defined in \eqref{d Poisson equation}. Here $\mathbf{x}=(x,y)^{T}$ and the forcing term $h(\mathbf{x})$ is given by 
\begin{equation}
\label{forcing_term}
h(\mathbf{x}) = \frac{2}{m}\sum_{i=1}^m (\omega_i \pi)^2 \sin(\omega_i\pi x) \sin(\omega_i \pi y),
\end{equation}
where $m$ is a positive integer. The corresponding exact solution to the problem ($\ref{d Poisson equation}$) is 
\begin{equation}
\label{Poisson_2d_sol}
u(\mathbf{x}) = \frac{1}{m}\sum_{i=1}^m \sin(\omega_i\pi x) \sin(\omega_i \pi y).
\end{equation}
In this example, we set $\omega_i = 2^i$ for $i = 1, 2, \dots, m$.

To test the performance of the modified HPKM-PINN with overlapping domain decomposition, we increase the complexity of the solution while refining the subdomains. Specifically, we fix the overlap ratio $\delta=2.9$ which is determined through testing to be well-suited for this problem. We then increment the parameter $m$ from 2 to 6 in the forcing term \eqref{forcing_term}, while adaptively increasing the number of subdomains and sampling points ($2^{m-1} \times 2^{m-1}$ subdomains and $(10 \times 2^{m-1}) \times (10 \times 2^{m-1})$ sampling points) to match the solution's growing complexity.

\begin{figure}[ht] 
    \centering
   \hspace{-0.1in}
   \begin{minipage}[b]{0.2 \textwidth}
   \includegraphics[width=\textwidth]{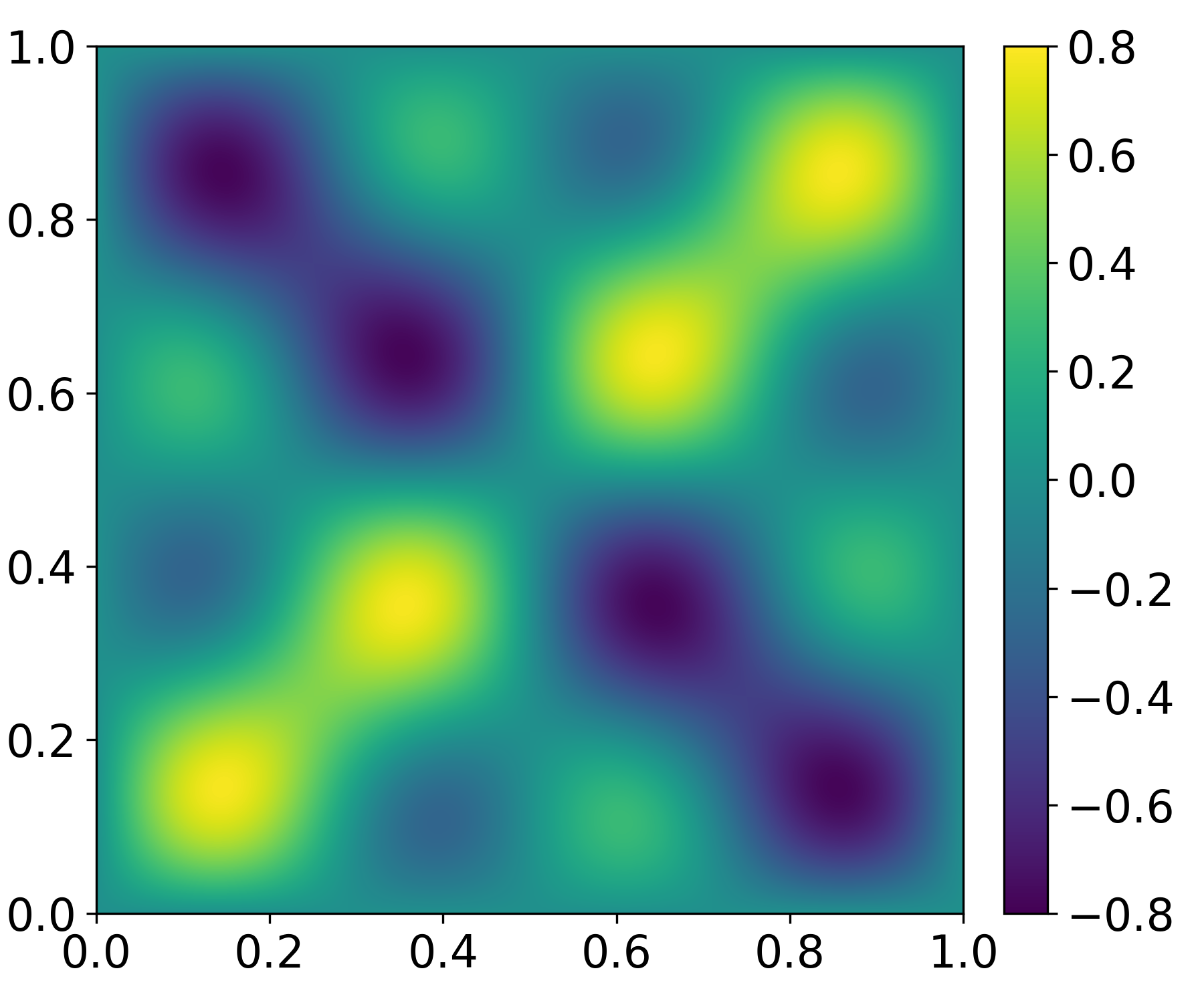}
   \begin{center}
    \hspace{-0.1cm}
   \footnotesize Exact solution \\ \footnotesize $m=2$
   \end{center}
   \end{minipage}
   \hspace{-0.1in}
   \begin{minipage}[b]{0.20 \textwidth}
   \includegraphics[width=\textwidth]{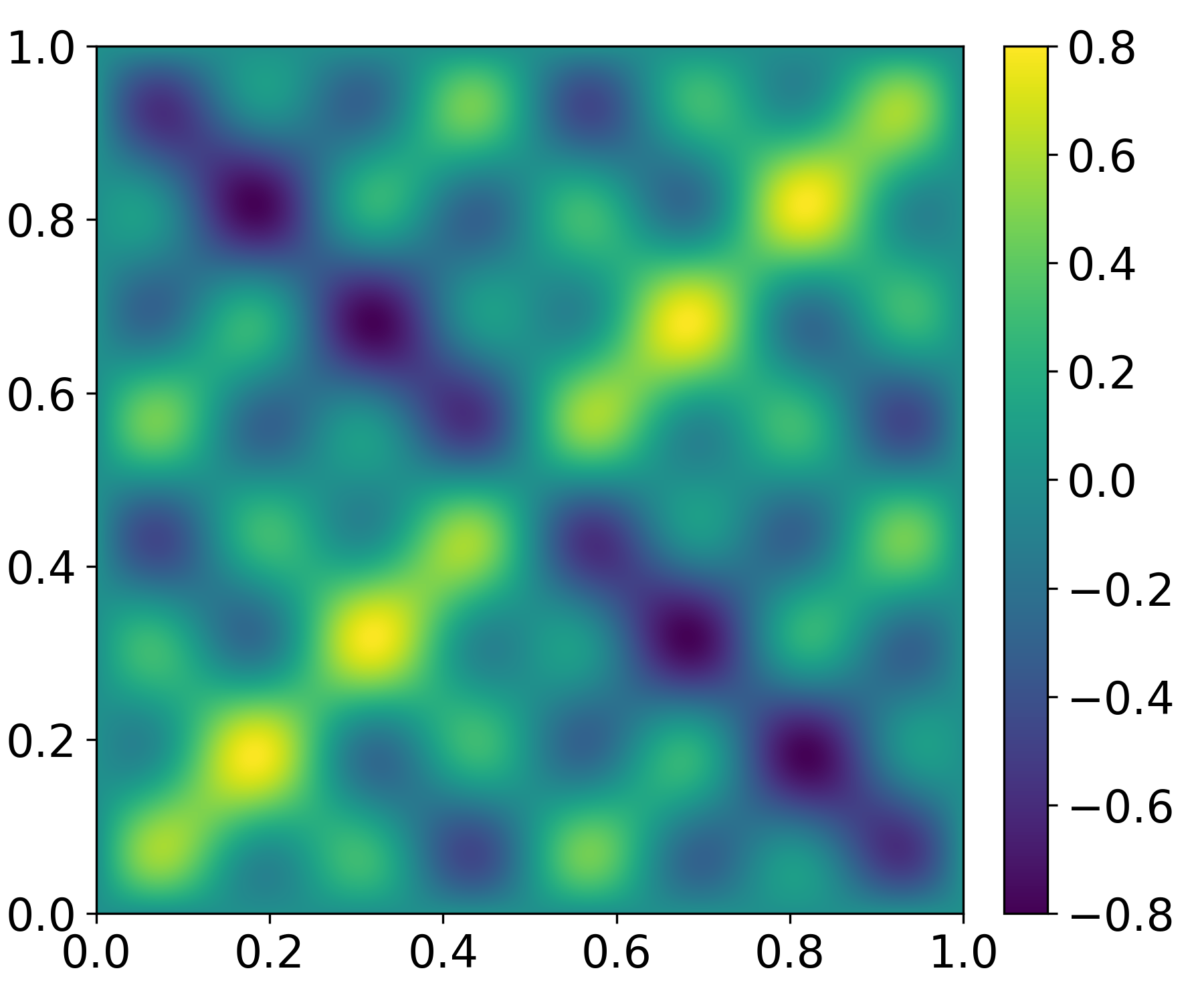}
   \begin{center}
   \footnotesize Exact solution \\ \footnotesize $m=3$
   \end{center}
   \end{minipage}
   \hspace{-0.1in}
   \begin{minipage}[b]{0.2 \textwidth}
   \includegraphics[width=\textwidth]{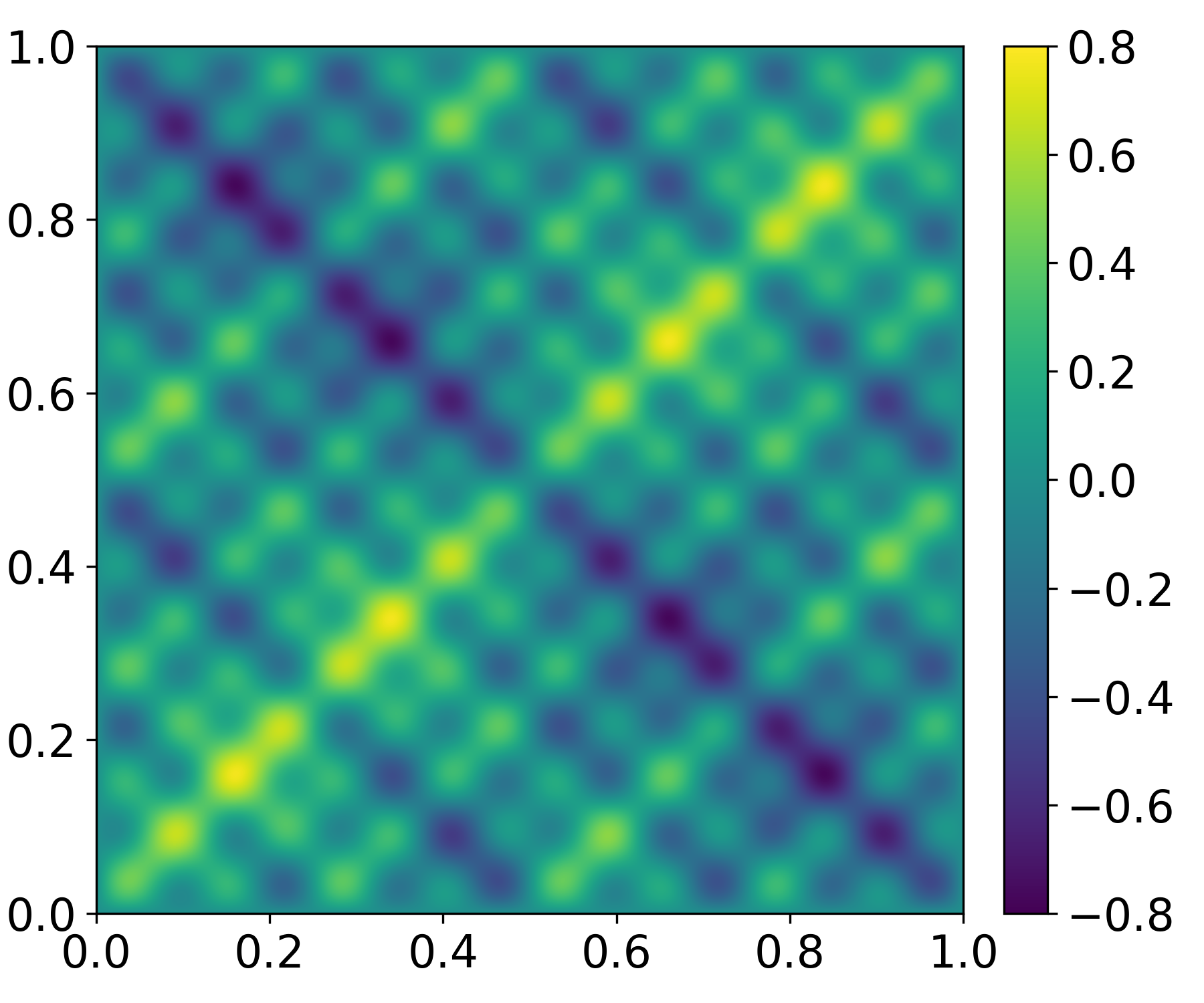}
   \begin{center}
   \footnotesize Exact solution \\ \footnotesize $m=4$
   \end{center}
   \end{minipage}
   \hspace{-0.1in}
   \begin{minipage}[b]{0.2 \textwidth}
   \includegraphics[width=\textwidth]{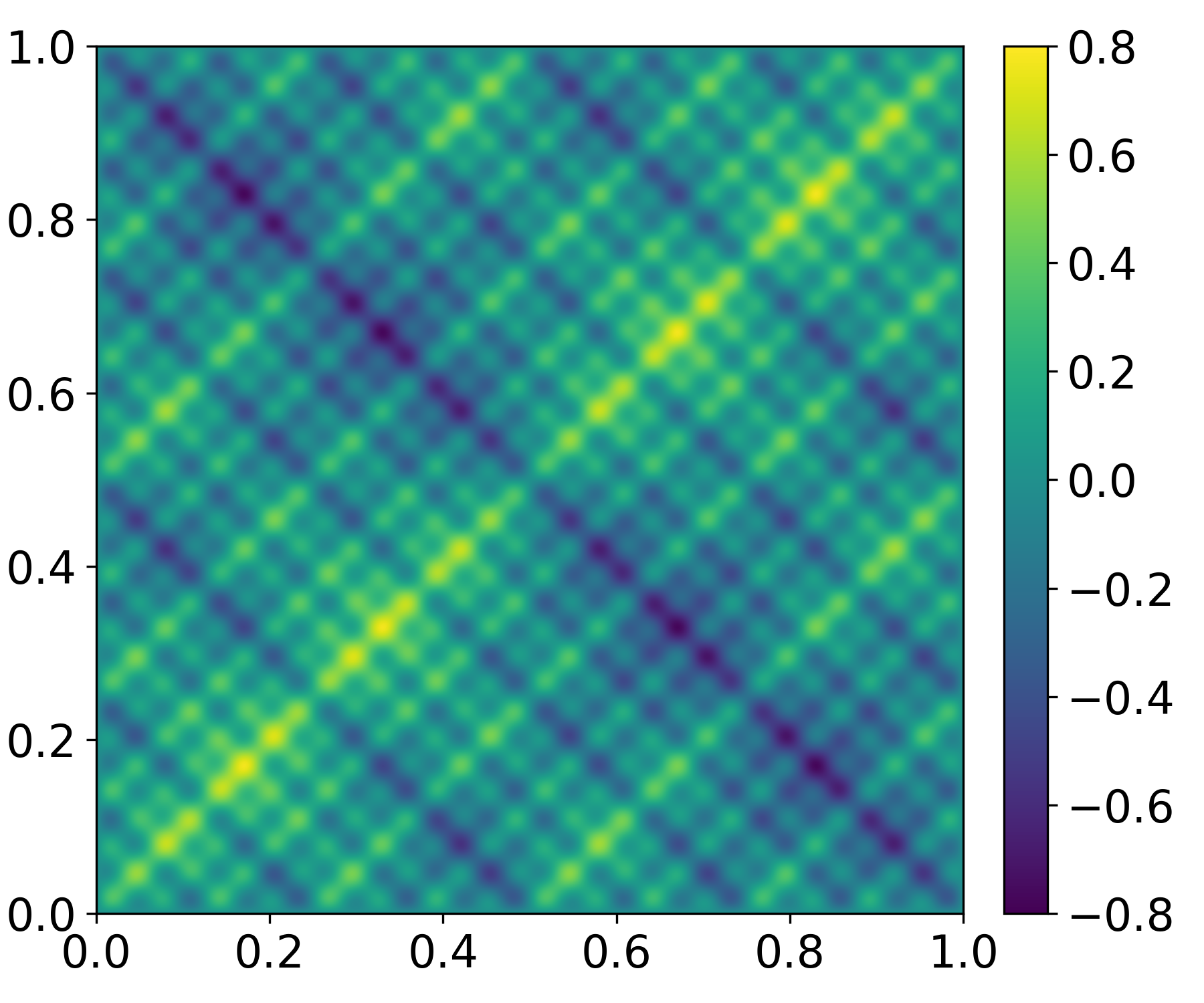}
   \begin{center}
   \footnotesize Exact solution \\ \footnotesize $m=5$
   \end{center}
   \end{minipage} 
 \hspace{-0.1in}
   \begin{minipage}[b]{0.20 \textwidth}
   \includegraphics[width=\textwidth]{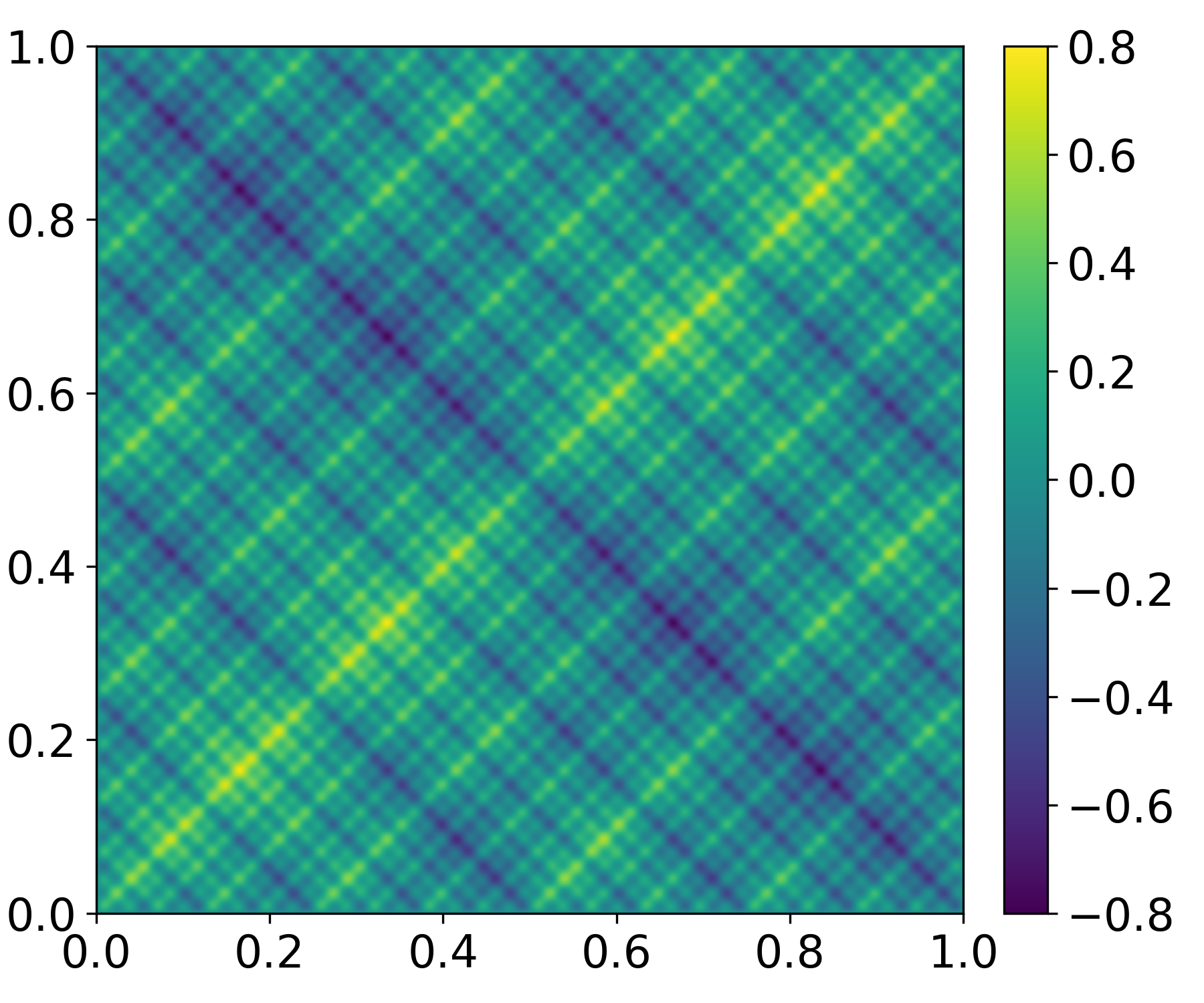}
   \begin{center}
   \footnotesize Exact solution \\ \footnotesize $m=6$
   \end{center}
   \end{minipage}
    \qquad
    \centering
   \hspace{-0.1in}
   \begin{minipage}[b]{0.2 \textwidth}
   \includegraphics[width=\textwidth]{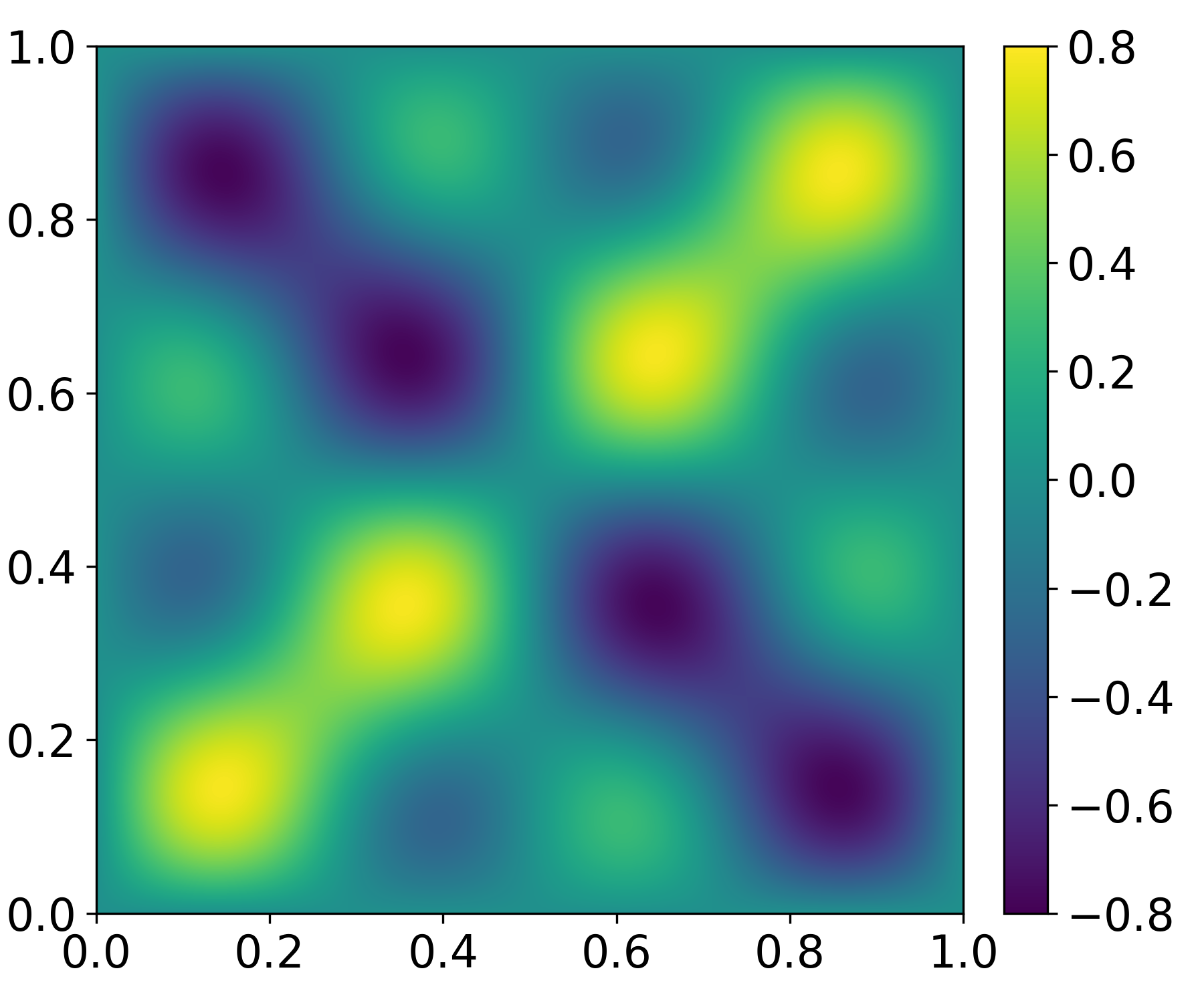}
   \begin{center}
    \footnotesize Predicted solution \\ \footnotesize Subdomain: $2\times 2$ \\ \footnotesize Point: $20 \times 20$ \\ $m = 2$
   \end{center}
   \end{minipage}
   \hspace{-0.1in}
   \begin{minipage}[b]{0.2 \textwidth}
   \includegraphics[width=\textwidth]{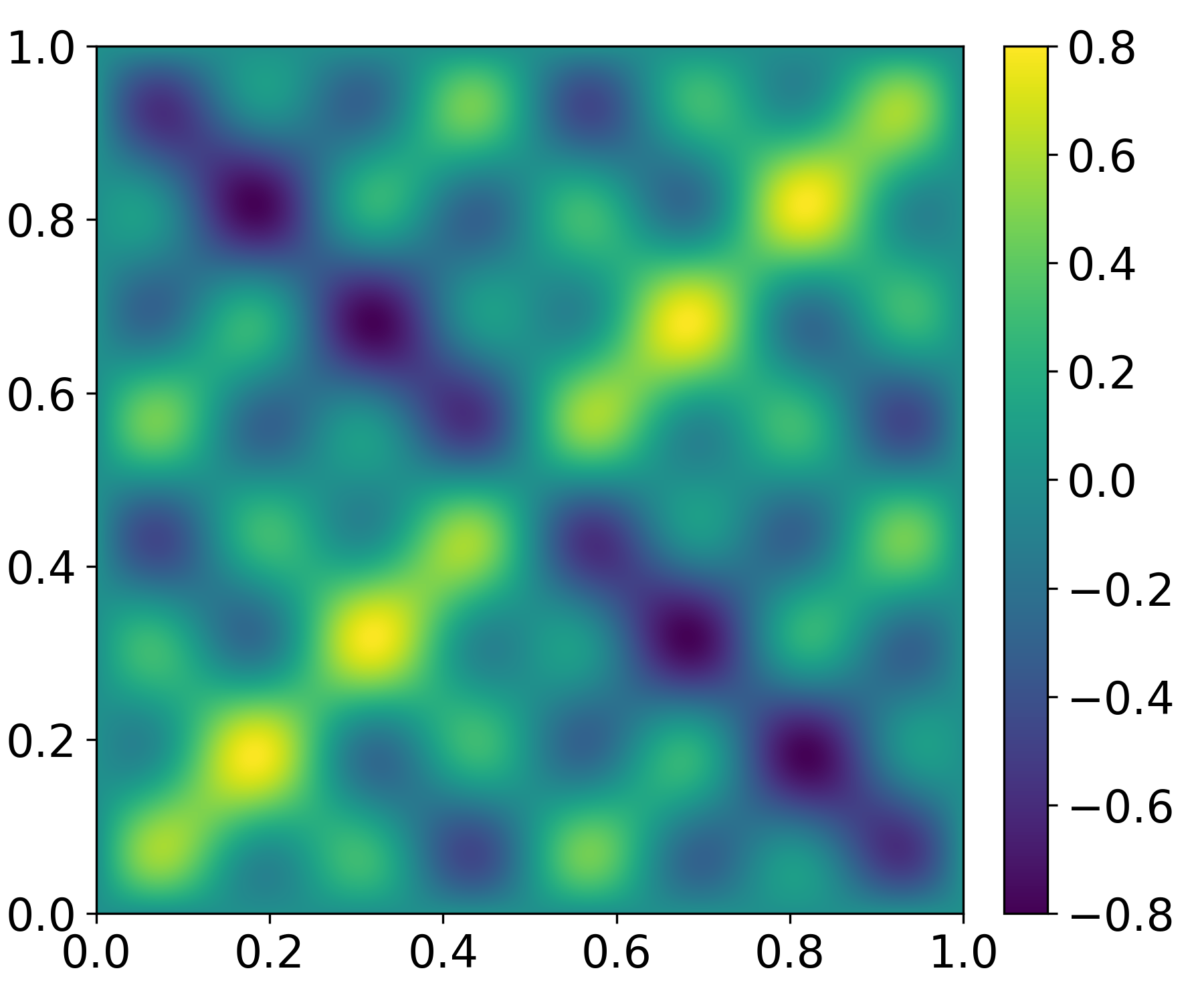}
   \begin{center}
    \footnotesize Predicted solution \\ \footnotesize Subdomain: $4\times 4$ \\ \footnotesize Point: $40  \times 40$ \\ $m = 3$
   \end{center}
   \end{minipage}
   \hspace{-0.1in}
   \begin{minipage}[b]{0.2 \textwidth}
   \includegraphics[width=\textwidth]{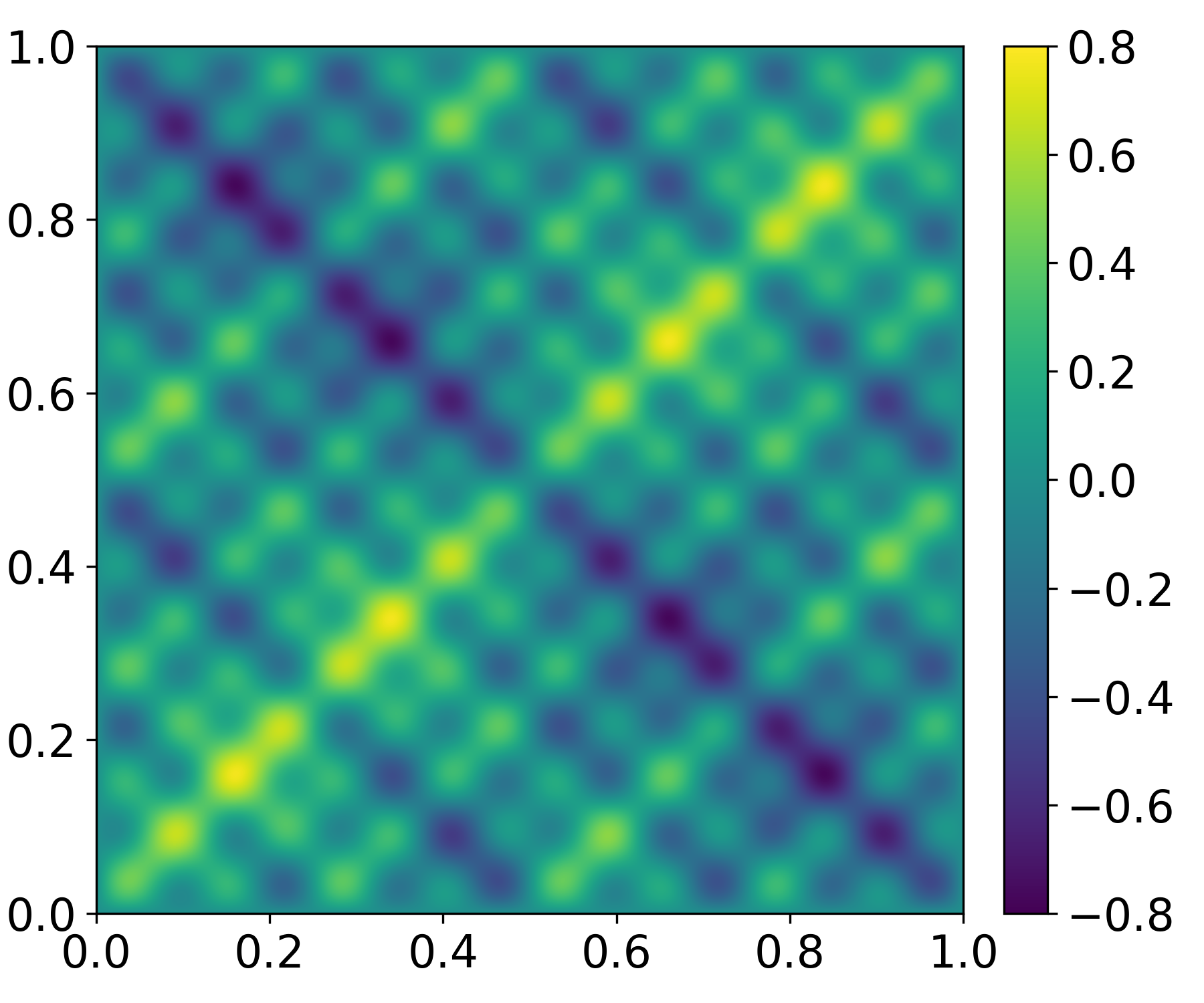}
   \begin{center}
    \footnotesize  Predicted solution \\ \footnotesize Subdomain: $8 \times 8 $\\ \footnotesize Point: $80 \times 80 $ \\ $m = 4$
   \end{center}
   \end{minipage}
   \hspace{-0.1in}
   \begin{minipage}[b]{0.2 \textwidth}
   \includegraphics[width=\textwidth]{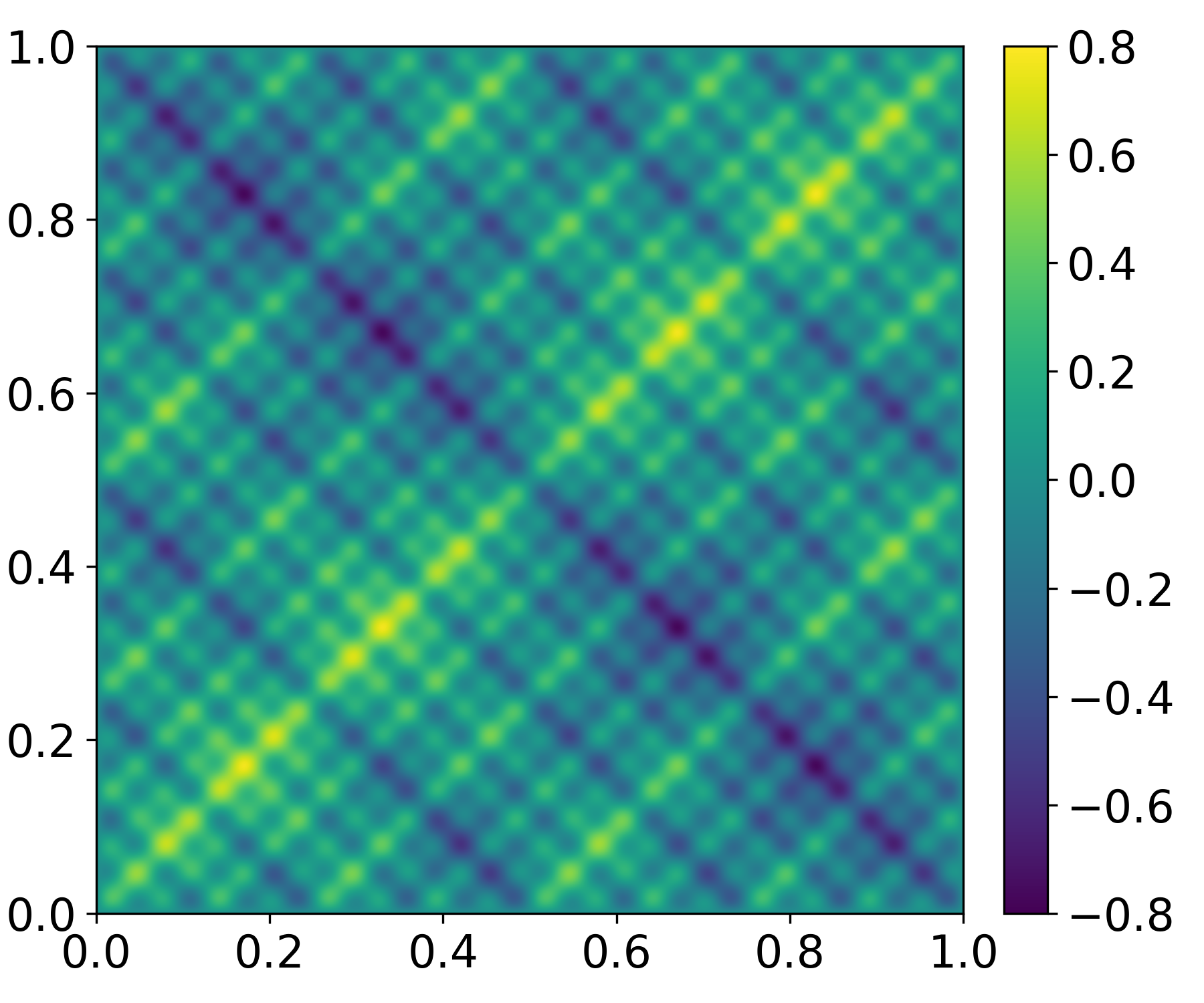}
   \begin{center}
    \footnotesize  Predicted solution \\ \footnotesize Subdomain: $16 \times 16$ \\ \footnotesize Point: $160 \times 160$ \\ $m = 5$
   \end{center}
   \end{minipage} 
 \hspace{-0.1in}
   \begin{minipage}[b]{0.2 \textwidth}
   \includegraphics[width=\textwidth]{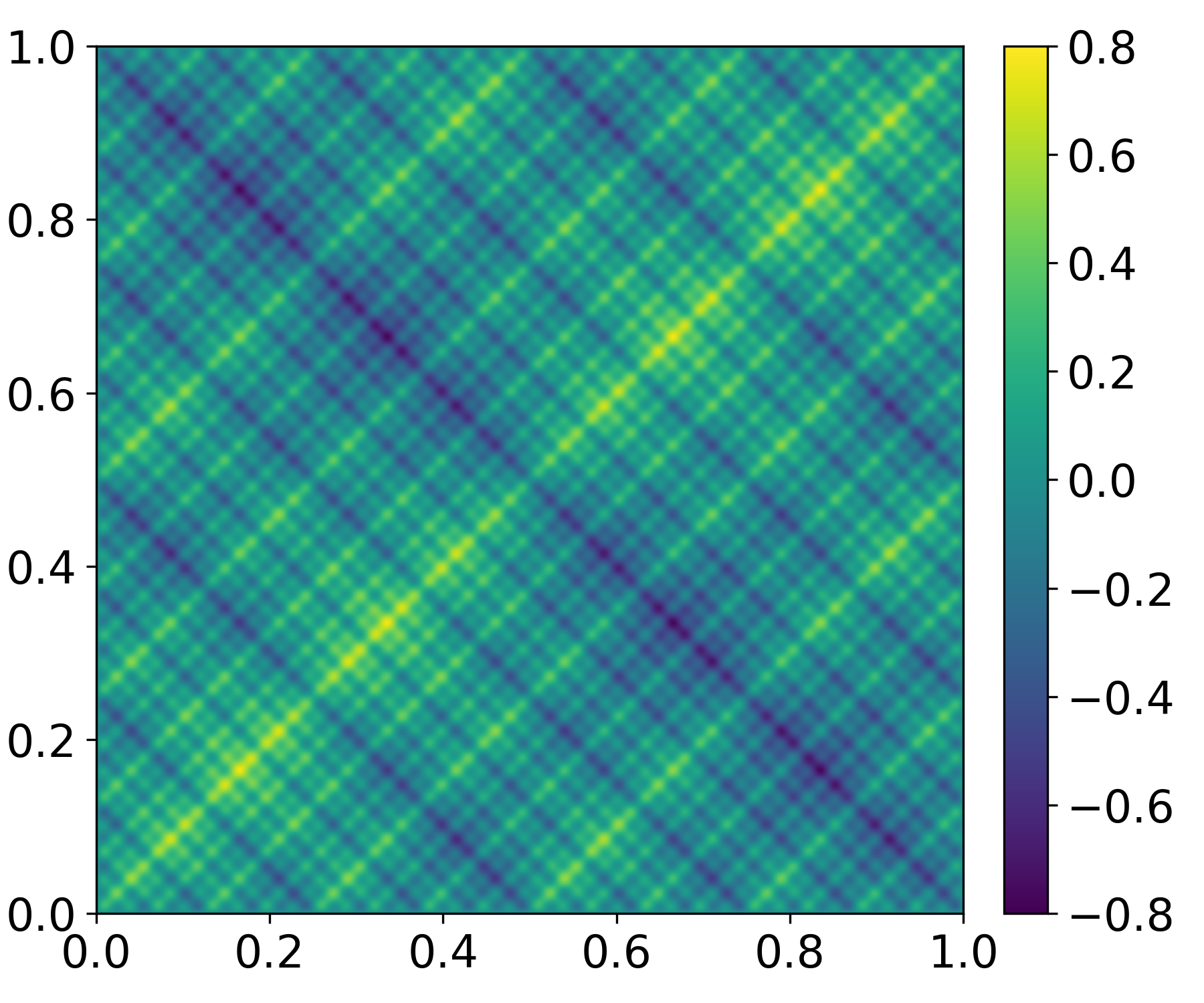}
   \begin{center}
    \footnotesize Predicted solution \\ \footnotesize Subdomain:$32 \times 32$  \\ \footnotesize Point: $320 \times 320$ \\ $m = 6$
   \end{center}
   \end{minipage}
 \qquad
    \centering
   \hspace{- 0.15 in}
   \begin{minipage}[b]{0.19 \textwidth}
   \includegraphics[width=1\textwidth]{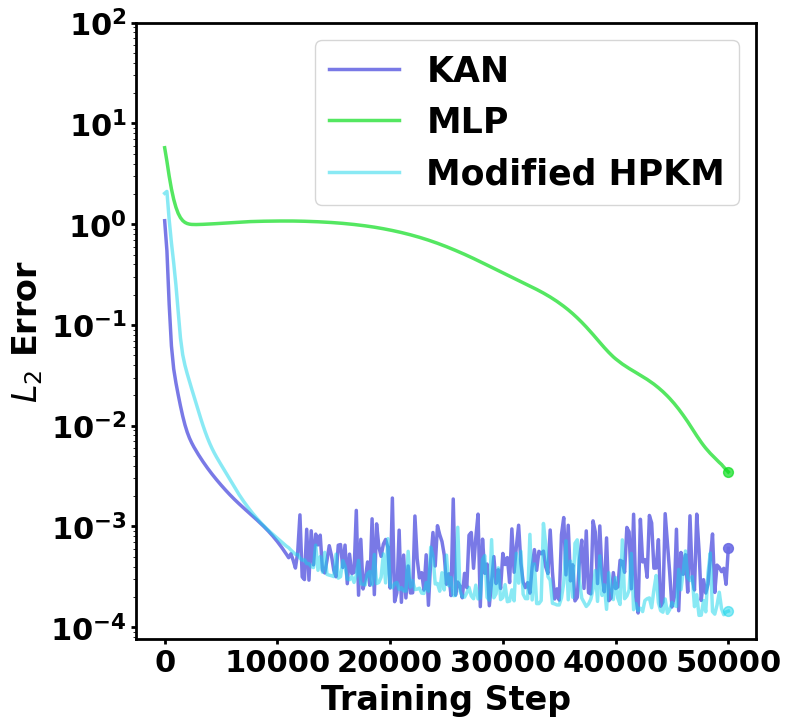}
   \begin{center}
   	\footnotesize $L_2$ error ($m = 2$)
   \end{center}
   \end{minipage}
   \hspace{-0.03  in}
   \begin{minipage}[b]{0.19\textwidth}
   \includegraphics[width=1\textwidth]{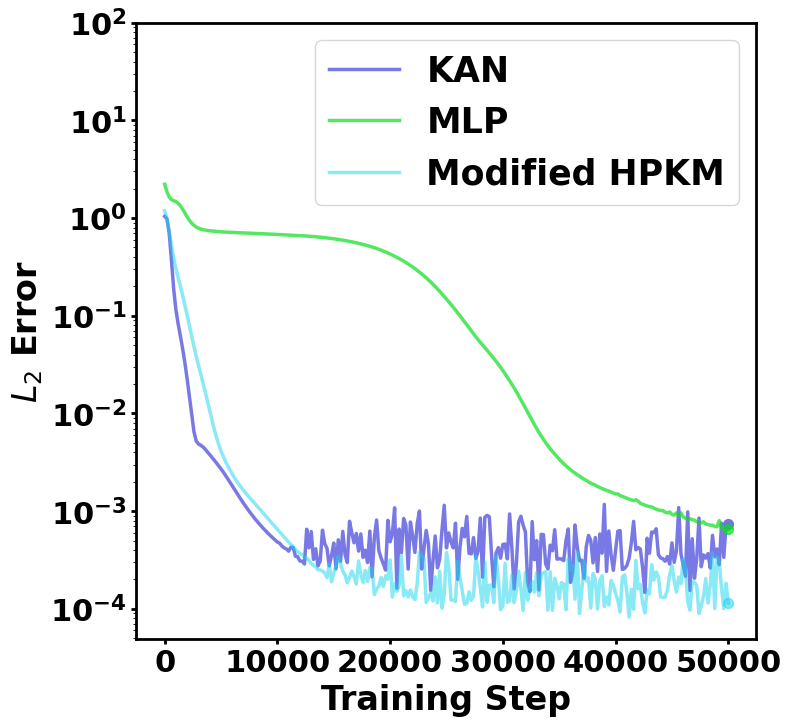}
    \begin{center}
   	\footnotesize $L_2$ error ($m = 3$)
   \end{center}
   \end{minipage}
   \hspace{-0.05 in}
   \begin{minipage}[b]{0.19 \textwidth}
   \includegraphics[width=1\textwidth]{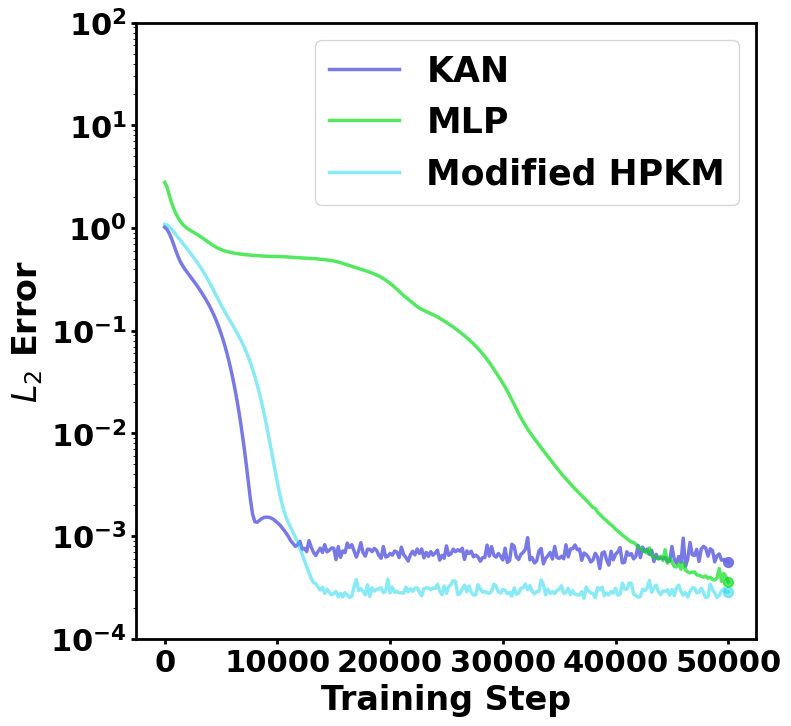}
    \begin{center}
   	\footnotesize $L_2$ error ($m = 4$)
   \end{center}
   \end{minipage}
   \hspace{-0.05 in}
   \begin{minipage}[b]{0.19 \textwidth}
   \includegraphics[width=1\textwidth]{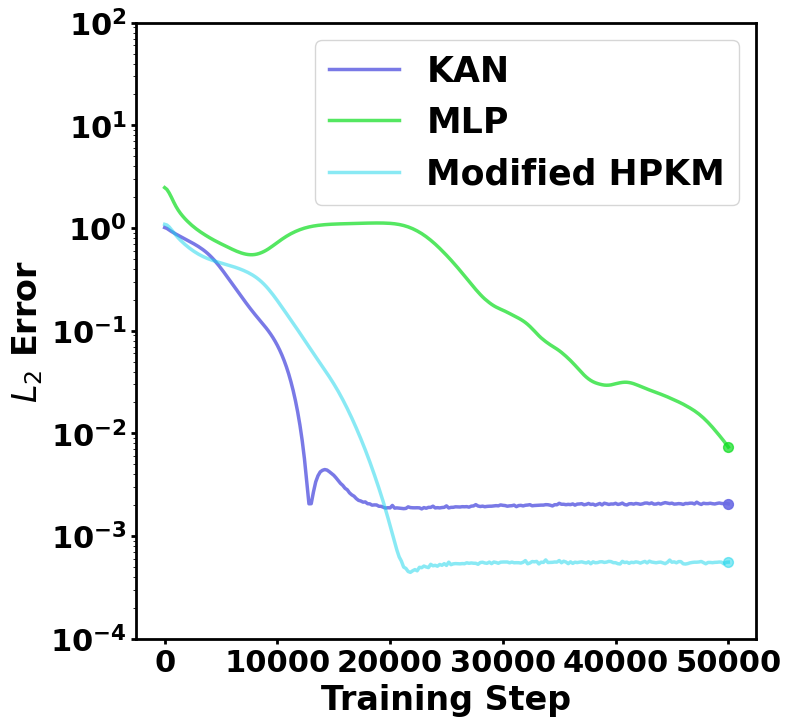}
    \begin{center}
   	\footnotesize $L_2$ error ($m = 5$)
   \end{center}
   \end{minipage} 
 \hspace{-0.05 in}
   \begin{minipage}[b]{0.19 \textwidth}
   \includegraphics[width=1\textwidth]{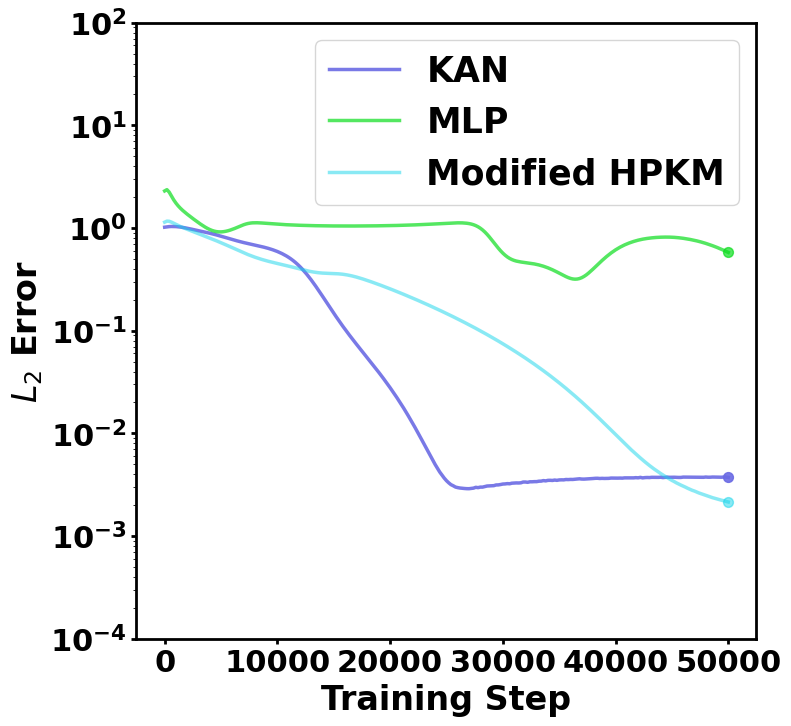}
    \begin{center}
   	\footnotesize $L_2$ error ($m = 6$)
   \end{center}
   \end{minipage}
   \begin{center}
   \caption{Two-dimensional Poisson equation with $m=2,3,4,5,6$. Top: Exact solutions. Middle: Predicted solutions obtained by the modified HPKM-PINN. Bottom: $L_2$ error comparison between KAN, MLP and modified HPKM architectures.}     \label{Hpkm-pinn-Poisson-equation}
   \end{center}
   \end{figure}

In Figure $\ref{Hpkm-pinn-Poisson-equation}$, the top and middle rows display the exact and predicted solution using the modified HPKM architecture, respectively. Then, we compare the $L_2$ error across multiple neural network models to analyze the influence of domain decomposition numbers and sampling points under varying parameter $m$ during network training, as shown in Figure $\ref{Hpkm-pinn-Poisson-equation}$. It can be observed in the bottom row of Figure $\ref{Hpkm-pinn-Poisson-equation}$ that both KAN and the modified HPKM rapidly converge to precision levels exceeding $10^{-3}$ compared to the MLP architecture. Furthermore, as $m$ increases, the problem complexity escalates and the $L_2$ error convergence results demonstrate enhanced stability. The final training errors presented in Table $\ref{err_total_result}$ demonstrate that the modified HPKM-PINN consistently achieves optimal accuracy. The results show that our hybrid architecture integrates the advantages of the MLP and KAN networks through adaptive weight adjustment, improving the convergence rate and the accuracy.

We remark that while our design adheres to the principles of $h$-$p$ refinement methods in the finite element method through a domain refinement strategy, a parameter discontinuity emerges when processing the solution \eqref{Poisson_2d_sol}. 
Specifically, as the parameter $m$ increases from $5$ to $6$, the corresponding number of subdomains in domain decomposition jumps from $256$ to $1024$. 
This discontinuity may introduce computational challenges, such as numerical instability and model overfitting. 
However, the modified HPKM-PINN with domain decomposition remains effective in solving the target problem under non-optimal parameters, achieving higher accuracy compared to baseline methods. 
The decrease of the $L_{2}$ error during training further validates the rationality of the architecture design and the generalization ability of our model.

\subsubsection{High-dimensional Poisson equation}

To further evaluate the performance of the modified HPKM-PINN with overlapping domain decomposition, we now consider a five-dimensional case ($d=5$) of the Poisson problem defined in \eqref{d Poisson equation}. The forcing term $h(\mathbf{x})$ in (\ref{d Poisson equation}) with $\mathbf{x}=(x_{1},x_{2},x_{3},x_{4},x_{5})^{T}$ is given by
\begin{equation*}
h(\mathbf{x}) = -5 \prod_{i=1}^5 \sin( \pi x_i),
\end{equation*}
and the corresponding exact solution is given by 
\begin{equation*}
u(\mathbf{x}) = -\frac{1}{\pi^2} \prod_{i=1}^5 \sin( \pi x_i).
\end{equation*}

\begin{figure}[htbp] 
	\centering
	\hspace{-0.1in}
	\begin{minipage}[b]{0.28\textwidth}
		\includegraphics[width=\textwidth]{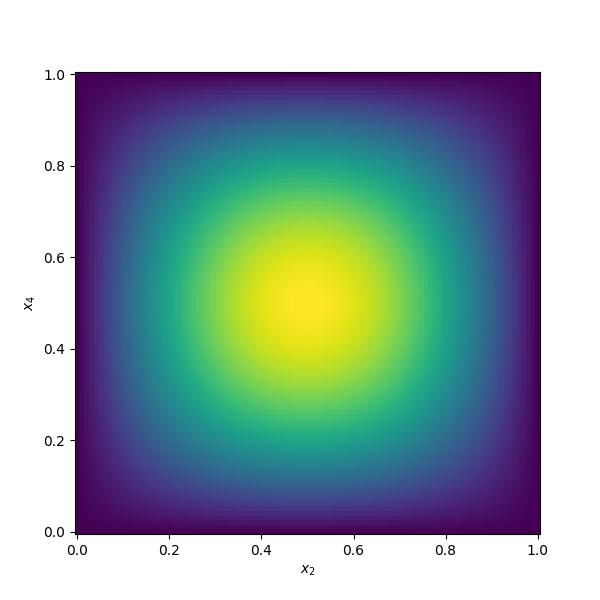}

        \subcaption{}
	\end{minipage}
	\hspace{-0.1in}
	\begin{minipage}[b]{0.28 \textwidth}
		\includegraphics[width=\textwidth]{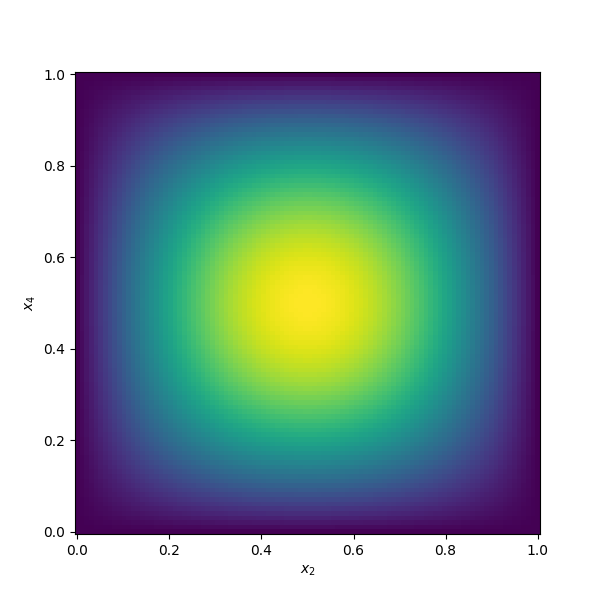}

	\subcaption{}
    \end{minipage}
	\qquad
	\centering
	\hspace{-0.1in}
	\begin{minipage}[b]{0.32 \textwidth}
		\includegraphics[width=\textwidth]{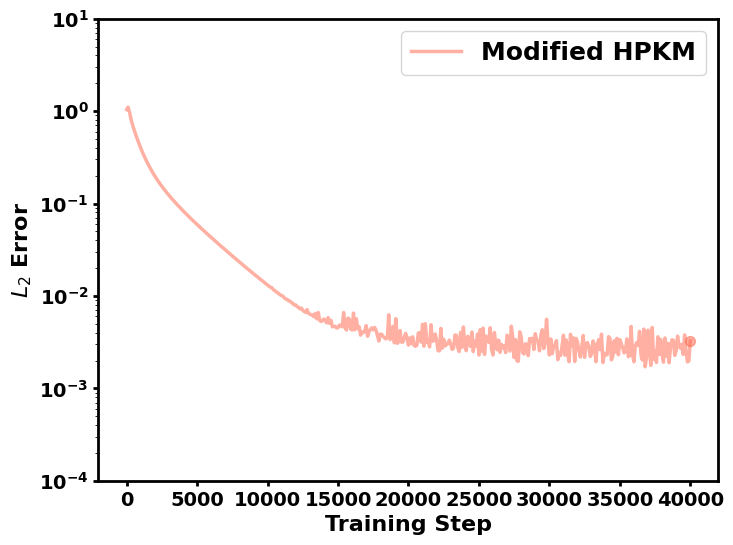} 

\subcaption{}
    \end{minipage}
	\begin{center}
		\caption{Five-dimensional Poisson equation: (a) A two-dimensional slice of the exact solution in the  $x_2$ and $x_4$ hyperplane; (b) A two-dimensional slice of the predicted solution obtained by the modified HPKM architecture in the  $x_2$ and $x_4$ hyperplane; (c) $L_2$ error by using the modified HPKM architecture.    }
		\label{5D-equation}
	\end{center}
\end{figure}

In solving the five-dimensional Poisson equation ($\ref{d Poisson equation}$), we partition the computational domain into $2\times2\times2\times2\times2$ subdomains with a sampling configuration of $10\times10\times10\times10\times10$ points and an overlap ratio of $\delta=1.9$. 
As shown in Figure $\ref{5D-equation}$, by using this simple domain decomposition, we obtain solutions of comparable accuracy. Notably, under the domain decomposition technique, higher-dimensional problems necessitate more subdomains, resulting in correspondingly higher memory requirements for computation. In conclusion, the modified HPKM architecture employs two single-layer network branches with reduced neuron counts. This design not only decreases GPU memory consumption but also helps solve the problem with minimal domain decomposition and sparse sampling configurations. The results in this part demonstrate the effectiveness of our method for solving high-dimensional problems.

\subsection{ Reaction-Diffusion equation}
The reaction-diffusion equation describes the evolution of chemical concentrations or population densities through coupled spatial diffusion and chemical interaction processes \cite{zengjia3-Reaction-diffusion}. 
The interaction of diffusion and reaction terms produces complex dynamics, such as traveling waves and stationary patterns. These behaviors are widely observed in fields such as ecology, chemistry, materials science, and neuroscience.
In this experiment, we consider the steady-state reaction-diffusion equation in one dimension for a scalar quantity, subject to homogeneous Dirichlet boundary conditions, that is,
\begin{equation}
	 \label{1-Reaction-diffusion}
\begin{aligned}
&D \Delta u(x) +\mu \tanh(u(x))=h(x),  ~~~&&x \in \Omega = (0,1),\\
&u(x)=0, &&x\in\partial \Omega, 
\end{aligned}
\end{equation}
where $u(x)$ represents the chemical concentration, $D$ denotes the diffusion coefficient and $\mu$ denotes the reaction rate. 
And $h(x)$ denotes the nonlinear reaction term, written as
\begin{equation*}
  h(x) = -3Dk^2 \pi^2( 2\sin(k\pi x)\cos^2(k\pi  x) - \sin^3(k\pi x) + \mu \tanh( \sin^3(k\pi x) ),
\end{equation*}
where $k$ is a positive constant. Thus, we have the chemical concentration 
\begin{equation*}
u(x) = \sin^3(k\pi x).                 
\end{equation*}
In this experiment, we set $D= 0.01$, $\mu = 0.7$, and $k = 8$.

\begin{figure}[htbp] 
    \centering
    \hspace{-0.1in}
   \begin{minipage}[b]{0.35\textwidth}
   \includegraphics[width=\textwidth]{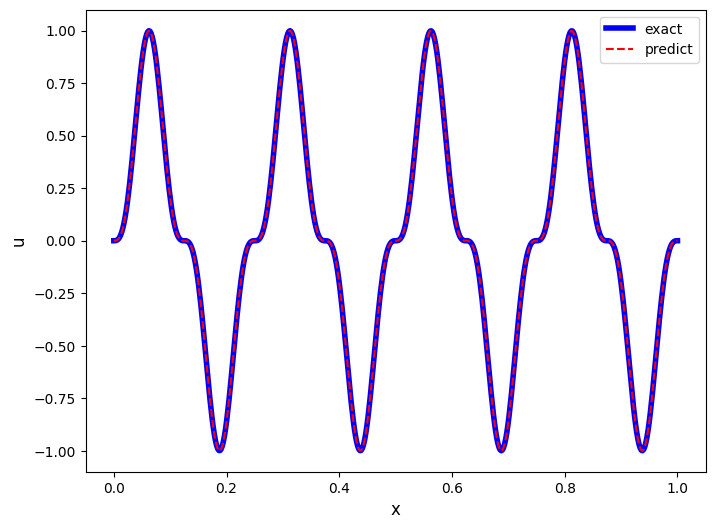} 
   \subcaption{}
   \end{minipage}
  \qquad
    \centering
   \hspace{0.3in}
   \begin{minipage}[b]{0.35 \textwidth}
   \includegraphics[width=\textwidth]{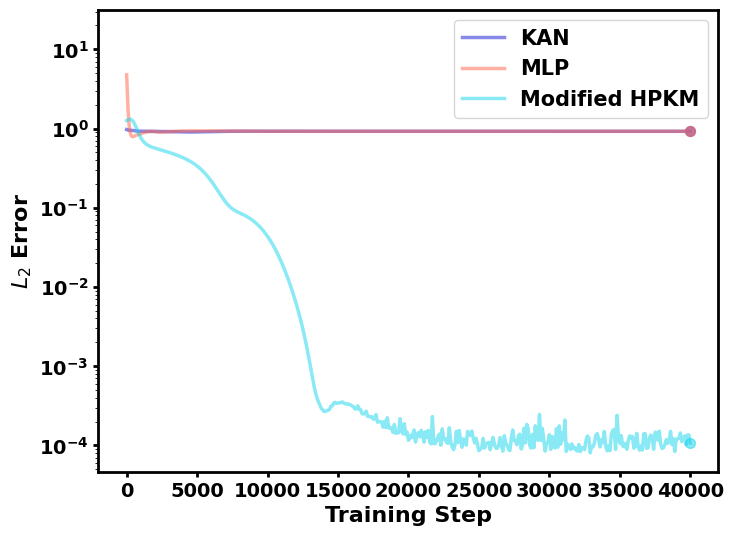} 
  \subcaption{}
   \end{minipage}
   \begin{center}
   \caption{Steady-State reaction-diffusion equation: (a) The exact and predicted solution obtained by the modified HPKM architecture; (b) $L_2$ error comparison between KAN, MLP and modified HPKM-PINN architectures.   }
   \label{Reaction-Diffusion equation text}
   \end{center}
   \end{figure}

In this numerical experiment, as shown in Figure $\ref{Reaction-Diffusion equation text}$ (a), the solution displays high-frequency oscillations, which presents a major challenge for standard PINNs in accurately capturing the solution features.
To solve this issue, we employ a configuration of 10 subdomains, 200 sampling points and the overlap ratio $\delta=2.9$. Through comparative analysis of three models based on overlapping domain decomposition techniques, we observe that both conventional MLP and KAN architectures with domain decomposition fail to accurately capture the high-frequency features of the solution in Figure $\ref{Reaction-Diffusion equation text}$ (b). In contrast, the modified HPKM-PINN demonstrates superior modeling accuracy. We remark that when increasing the number of subdomains, the MLP or KAN architecture can also achieve computational accuracy. This phenomenon further verifies the advantage of our proposed model in solving complex high-frequency problems while maintaining fewer parameters and higher accuracy.

In Table \ref{err_total_result}, the $L_{2}$ errors of the three network architectures after completing the training are presented, yielding 0.9226 (MLP), 0.9225 (KAN), and $1.07\times 10^{-4}$ (modified HPKM). These results clearly demonstrate that the modified HPKM-PINN, combined with overlapping domain decomposition, achieves higher accuracy than standalone MLP or KAN architectures. It confirms that our method, as a novel architecture, provides a framework for solving partial differential equations using neural networks combined with domain decomposition. For future work, we will investigate the applicability of this architecture to other types of partial differential equations.

\subsection{Allen--Cahn equation} 
Allen--Cahn equation is a phase-field model proposed by Allen and Cahn to describe interfacial evolution dynamics during phase separation processes \cite{zengjia4Allen--Cahn}. This equation plays a significant role in modeling incompressible, immiscible two-phase fluid systems and has been used to investigate spinodal decomposition with hydrodynamic effects as well as the dynamic behavior of polymeric fluids in multiphase systems.
We consider the following Allen--Cahn equation,
\begin{equation}
\label{eq:ac}
\begin{aligned}
&\frac{\partial u}{\partial t} = \epsilon^2 \frac{\partial^2 u}{\partial x^2} + u - u^3 + h(x,t),  ~~~&&x \in (0,1), ~~ t \in (0,1],\\
& u(0,t) = u(1,t) = 0,   \\
& u(x,0) = \sin(\pi x) + 0.1\sin(10\pi x), 
\end{aligned}
\end{equation}
where $u = u(x, t)$ is a phase variable and $\epsilon$ indicates the thickness of the interface. $h(x, y)$ is a forcing term expressed as 
\begin{equation*}
\begin{split}
h(x,t) = & -2\pi \sin(\pi x)\sin(2\pi t) - 2\pi \sin(10\pi x)\sin(20\pi t) \\
         & + \epsilon^2 \left( \pi^2 \sin(\pi x)\cos(2\pi t) + 10\pi^2 \sin(10\pi x)\cos(20\pi t) \right) \\
         & + \left( \sin(\pi x)\cos(2\pi t) + 0.1\sin(10\pi x)\cos(20\pi t) \right)^3 \\
         & - \left( \sin(\pi x)\cos(2\pi t) + 0.1\sin(10\pi x)\cos(20\pi t) \right).
\end{split} 
\end{equation*}
Then the exact solution $u$ to the Allen--Cahn equation \eqref{eq:ac} is given in the following form 
\begin{align*}
u(x,y) =  \sin(\pi x) \cos(2\pi t) + 0.1 \sin(10\pi x) \cos(20\pi t).
\end{align*}
Here we set $\epsilon = 0.1$.

 \begin{figure}[htbp] 
    \centering
   \hspace{-0.1in}
   \begin{minipage}[b]{0.32\textwidth}
   \includegraphics[width=\textwidth]{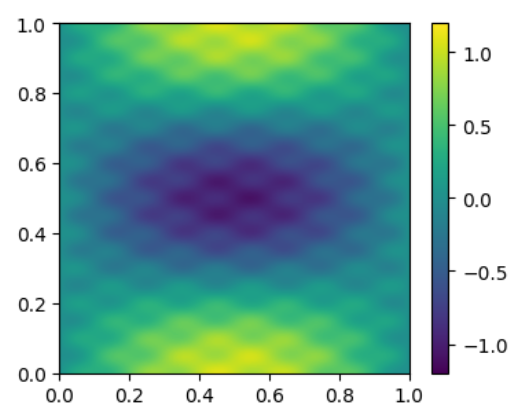}
   \subcaption{}
   \end{minipage}
   \hspace{-0.1in}
   \begin{minipage}[b]{0.32 \textwidth}
   \includegraphics[width=\textwidth]{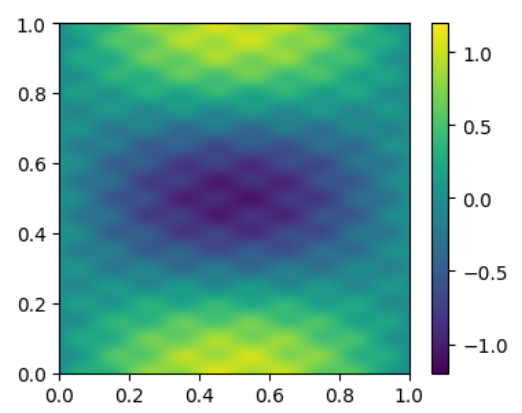}
   \subcaption{}
   \end{minipage}
  \qquad
    \centering
   \hspace{-0.1in}
   \begin{minipage}[b]{0.32 \textwidth}
   \includegraphics[width=\textwidth]{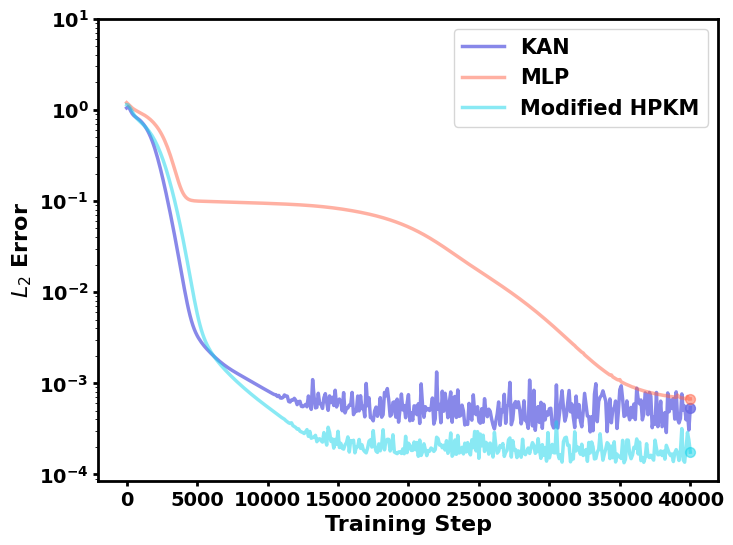} 
   \subcaption{}
   \end{minipage}
   \begin{center}
   \caption{Allen--Cahn equation: (a) Exact solution; (b) Predicted solution obtained by the modified HPKM-PINN; (c) $L_2$ error comparison between KAN, MLP and modified HPKM architecture.}
    \label{AC-equation}
   \end{center}
   \end{figure}

In this experiment, we incorporate the temporal dimension into the domain decomposition framework by employing a ($260 \times 260$) grid of sampling points, partitioned into $20 \times 20$ subdomains with an overlap ratio of $\delta = 2.9$. The results in Figure $\ref{AC-equation}$ (a)-(b) confirm that our parameter configuration provides an accurate prediction of the solution for the Allen--Cahn equation \eqref{eq:ac}. Detailed error convergence patterns in Figure $\ref{AC-equation}$ (c) reveal three aspects of the performance of the three models: (i) all the three architectures (MLP, KAN, and modified HPKM) with overlapping domain decomposition reach comparable solution accuracy, confirming their ability to solve this type of equation; (ii) the models using either KAN or modified HPKM architectures demonstrate faster $L_2$ error decay during the initial training phase, indicating their optimization dynamics during early learning stages; (iii) The modified HPKM architecture achieves the lowest error among all tested models in solving the Allen--Cahn equation \eqref{eq:ac}. As is shown in Table \ref{err_total_result} and Figure $\ref{AC-equation}$, compared with MLP and KAN architecture, the modified HPKM architecture achieves a higher numerical accuracy, providing a robust solution methodology for partial differential equations characterized by mixed high-frequency and low-frequency features. The superior performance of the modified HPKM architecture highlights the advantages of hybrid architectural design for complex flow modeling applications.

\begin{table}[h] 
\caption{$L_{2}$ error comparison across different network architectures.}
\label{err_total_result}
\renewcommand{\arraystretch}{2}  
\begin{center}
  \begin{tabular} 
               {@{} c  c |   c| c |c @{}} 
    \hline           
 \multirow{2}{*}{PDE}&  & \multicolumn{3}{c}{\textbf{Architecture}} \\ \cline{3-5}
&  & \textbf{MLP} & \textbf{ KAN} & \textbf{Modified HPKM} \\
    \hhline{|-----|}
    \multirow{3}{*}{Helmholtz (\ref{Hemholtz equation}) } & $\omega=16$     & $2.161 \times 10^{-3}$ & $5.115 \times 10^{-4}$ & $\mathbf{2.362 \times 10^{-4}}$ \\
    \hhline{|~|----|}
    & $\omega=32$     & $1.134 \times 10^{-1}$ & $6.932 \times 10^{-4}$ & $\mathbf{3.372 \times 10^{-4}}$ \\
    \hhline{|~|----|}
    & $\omega=48$     & $2.714 $ & $1.124 \times 10^{-3}$ & $\mathbf{8.698 \times 10^{-4}}$ \\
    \hhline{|-----|}
    \multirow{6}{*}{2D Poisson (\ref{d Poisson equation})}   & $m=2$     & $3.469 \times 10^{-3}$  & $6.066 \times 10^{-4}$ & $\mathbf{1.442 \times 10^{-4}}$ \\
    \hhline{|~|----|}
                               & $m=3$     & $6.510 \times 10^{-4}$  & $7.411\times 10^{-4}$ & $\mathbf{1.134 \times 10^{-4}}$ \\
    \hhline{|~|----|}
                             & $m=4$    & $3.586 \times 10^{-4}$ & $5.568 \times 10^{-4}$ & $\mathbf{2.868 \times 10^{-4}}$ \\
    \hhline{|~|----|}
                            & $m=5$    & $7.361 \times 10^{-3}$ & $2.058 \times 10^{-3}$ & $\mathbf{5.571 \times 10^{-4}}$ \\
    \hhline{|~|----|}
                             & $m=6$      & $5.787 \times 10^{-1}$ & $3.72 \times 10^{-3}$ & $\mathbf{2.147 \times 10^{-3}}$ \\
    \hhline{|-----|}
   Reaction diffusion (\ref{1-Reaction-diffusion})  &      & $9.226 \times 10^{-1}$ & $9.225 \times 10^{-1}$ & $\mathbf{1.070 \times 10^{-4} }$\\
    \hhline{|-----|}
    Allen--Cahn (\ref{eq:ac})  &       & $6.722 \times 10^{-4}$ & $5.397 \times 10^{-4}$ & $\mathbf{1.751 \times 10^{-4}}$ \\
    \hhline{|-----|}
  \end{tabular}
     \vspace{0.2cm} 

\end{center}
 \end{table}

\section{Conclusion}

\textcolor{black}{This study proposes a modified HPKM-PINN that incorporates overlapping domain decomposition via partition of unity functions. Within this framework,} we construct a convex combination of MLP and KAN outputs using trainable weight parameters optimized by S-shaped functions to balance their contributions in handling different frequency features and reduce computational costs. By employing overlapping domain decomposition, the global problem is divided into parallel subproblems, significantly improving convergence efficiency and alleviating the challenges of global optimization, providing an efficient framework for solving multi-scale and high-frequency partial differential equations. 

The results of extensive numerical tests in this work demonstrate the superior performance of the modified HPKM-PINN integrated with overlapping domain decomposition. The modified HPKM-PINN achieves higher accuracy than single MLP or KAN architecture. Although the model requires higher training costs than traditional MLP architectures due to its larger parameter size, it achieves a substantial parameter reduction compared to single KAN implementations.

In conclusion, we propose a modified HPKM-PINN computational framework that incorporates overlapping domain decomposition techniques, establishing an effective computational paradigm for solving complex high-frequency multiscale problems. This framework exhibits strong extensibility, enabling integration with residual attention weighting, multi-level FBPINNs, and Fourier feature embeddings for improved training. Future directions include evolutionary pattern-based optimization of subdomain parameters.

\section{Acknowledgements}
Q. Huang was partially supported by the  National Natural Science Foundation of China (No.12371385), and Y. Zhao was partially supported by the National Natural Science Foundation of China (No.12401509).

\bibliographystyle{plain}

\end{document}